\documentclass[aap]{imsart}

\RequirePackage{amsthm,amsmath,amsfonts,amssymb}
\RequirePackage[authoryear]{natbib}
\RequirePackage[colorlinks,citecolor=blue,urlcolor=blue]{hyperref}
\usepackage{comment,enumitem}
\startlocaldefs
\theoremstyle{plain}
\newtheorem{theorem}{Theorem}
\newtheorem{lemma}{Lemma}
\newtheorem{proposition}{Proposition}
\newtheorem{corollary}{Corollary}
\theoremstyle{definition}
\newtheorem{remark}{Remark}
\newtheorem*{condition}{Condition}
\newtheorem{assumption}{Assumption}
\newcommand\norm[1]{\left\lVert#1\right\rVert}
\newcommand\abs[1]{\left\lvert#1\right\rvert}
\newcommand\prs[1]{\left(#1\right)}
\newcommand\sbk[1]{\left[#1\right]}

\newcommand{\defeq}{\stackrel{\Delta}{=}}
\newcommand{\E}{\mathbb{E}}

\newcommand{\X}{\mathcal{X}}
\newcommand{\V}{\mathcal{V}}

\newcommand{\R}{\mathbb{R}}

\newcommand{\gti}{\rightarrow\infty}
\newcommand{\gtz}{\rightarrow0}

\endlocaldefs

\begin{document}

\begin{frontmatter}
\title{Computable Bounds on Convergence of Markov Chains in Wasserstein Distance via Contractive Drift}
\runtitle{Wasserstein Convergence Bounds for Markov Chains}

\begin{aug}
\author{\fnms{Yanlin}~\snm{Qu}\ead[label=e1]{quyanlin@stanford.edu}},
\author{\fnms{Jose}~\snm{Blanchet}\ead[label=e2]{jose.blanchet@stanford.edu}}
\and
\author{\fnms{Peter}~\snm{Glynn}\ead[label=e3]{glynn@stanford.edu}}
\address{Department of Management Science and Engineering,
Stanford University\printead[presep={,\ \\}]{e1,e2,e3}}

\end{aug}

\begin{abstract}
We introduce a unified framework to estimate the convergence of Markov chains to equilibrium in Wasserstein distance. The framework can provide convergence bounds with rates ranging from polynomial to exponential, all derived from a {\it contractive drift} condition that integrates not only contraction and drift but also coupling and metric design. The resulting bounds are computable, as they contain simple constants, one-step transition expectations, but no equilibrium-related quantities.
We introduce the large M technique and the boundary removal technique to enhance the applicability of the framework, which is further enhanced by deep learning in \cite{qu2024deep}. 
We apply the framework to non-contractive or even expansive Markov chains arising from queueing theory, stochastic optimization, and Markov chain Monte Carlo.
\end{abstract}



\end{frontmatter}

\section{Introduction}
 The long-term equilibrium of general state-space Markov chains is crucial in a wide array of applications. It is particularly relevant to inference and 
 stochastic optimization algorithms, such as Markov chain Monte Carlo and 
 constant step-size stochastic gradient descent. Additionally, it plays a 
 significant role in diverse stochastic models utilized across engineering, and the social and physical sciences, encompassing areas like logistics, supply 
 chains, economics, and population dynamics.

 The total variation (TV) distance has long been the standard metric used to measure the convergence of Markov chains. Significant efforts have been made to establish TV convergence bounds (\cite{meyn1994computable,tuominen1994subgeometric,rosenthal1995minorization,jarner2002polynomial,douc2004practical,baxendale2005renewal,hairer2011yet,andrieu2015quantitative,zhou2022dimension}). In the context of general state-space Markov chains, these types of bounds typically involve verifying so-called drift and minorization (D\&M) conditions; see \cite{meyn_tweedie_glynn_2009}. The drift condition ensures that the Markov chain moves towards a selected region. On such a region, the minorization condition guarantees a mixture representation for the transition kernel which can be used to obtain a suitable coupling with a stationary version of the chain. This coupling analysis, which is essential for TV bounds, tends to produce estimates that may be too conservative for practical use in applications such as those mentioned earlier (\cite{jones2001honest}), especially in high-dimensional settings (\cite{rajaratnam2015mcmc,qin2021limitations}). 
 
 An alternative to the TV distance is the so-called Wasserstein distance (\cite{villani2009optimal}). The Wasserstein distance measures the minimal expected cost (minimizing over joint distributions preserving source and target marginals) of transporting mass encoded by one distribution (the source) into another distribution (the target). The cost per unit of transported mass from a source location to a target location is measured using a distance between locations (typically a norm in a Euclidean space). It is well known that under mild integrability conditions, the Wasserstein distance metrizes the weak convergence topology and therefore it is weaker than the TV distance. However, as we shall see, in many applications of interest, the Wasserstein distance convergence provides a more effective tool for studying rates of convergence to stationarity. The quality of the estimates also seems to improve compared to those obtained via TV, when applying to high-dimensional models (\cite{hairer2011asymptotic,hairer2014spectral,durmus2015quantitative,mangoubi2019mixing,qin2022wasserstein}). In fact, an interesting and topical example in which convergence may fail in TV (as we shall discuss in Section \ref{section_tvw}) arises in the context of constant step-size stochastic gradient descent. Simply put, the Wasserstein distance provides a criterion that is powerful enough to quantify convergence, yet versatile enough to be applicable to practical chains that may take a very long time to converge in TV or that may not even converge at all in TV. Because of these benefits, the Wasserstein distance as a measure of convergence to equilibrium has steadily gained popularity over the years (\cite{gibbs2004convergence,ollivier2009ricci,madras2010quantitative,hairer2011asymptotic,butkovsky2014subgeometric,durmus2015quantitative,durmus2016subgeometric,douc2018markov,biswas2019estimating,eberle2019quantitative,butkovsky2020generalized,qin2022geometric,qin2022wasserstein,sandric2022subexponential}). These methods typically involve replacing the minorization condition (D\&M: drift and minorization) with a contraction condition (D\&C: drift and contraction); see, e.g., \cite{hairer2011asymptotic}.
 Establishing a contraction condition sometimes requires designing a new coupling or a new metric; see, e.g., \cite{butkovsky2020generalized}.
 In general, applying these methods to quantitatively analyze the convergence of realistic Markov chains is challenging, but for Langevin algorithms and Hamiltonian Monte Carlo, quantitative results have been obtained where both couplings and metrics are carefully designed to establish contraction (\cite{durmus2015quantitative,mangoubi2019mixing,eberle2019couplings,bou2020coupling,monmarche2021high}). 
 There is also a parallel line of research on bounding the convergence of diffusions in Wasserstein distance (\cite{hairer2008spectral,eberle2011reflection,eberle2016reflection,zimmer2017explicit,eberle2019sticky,eberle2019quantitativeH,lazi2021sub,nguyen2024polynomial}). 
 
 We briefly mention a different (but very powerful) set of methods involving stochastic localization and spectral independence techniques for bounding mixing times of discrete Markov chains; see \cite{Ronen22}, \cite{Anari20}. We leave the connections between these methods (typically designed for total variation convergence) and the methods proposed in this paper (utilizing intrinsic metrics and Wasserstein distance) as an interesting topic for future research.

\subsection{Main contributions}
 The objective of this paper is to introduce a unified framework to quantitatively bound the convergence of Markov chains in Wasserstein distance. This framework integrates drift, contraction, coupling and metric design into a single inequality, termed the {\it contractive drift} (CD) condition. We devise several techniques to establish CDs for a wide range of (not necessarily contractive) examples in queueing theory, stochastic optimization, and Markov chain Monte Carlo. For these examples, we obtain sharp or even parametrically sharp convergence bounds. More importantly, this CD framework serves as the theoretical foundation of the {\it Deep Contractive Drift Calculator} (DCDC), the first general-purpose sample-based algorithm to bound the convergence of Markov chains, which is recently introduced in \cite{qu2024deep}. As its name suggests, DCDC is powered by deep learning, which brings the convergence analysis of Markov chains from the pen-and-paper age to the era of AI. In the current paper, we focus on developing the CD framework and illustrating its effectiveness as an analytical framework. Specifically, our primary contributions include the following:

\begin{enumerate}[label=\roman*)]
    \item We introduce the notion of ``contractive drift'' (CD) that we utilize to derive explicit convergence bounds for Markov chains that exhibit polynomial (Theorem \ref{theorem_polynomial_b}), semi-exponential (Theorem \ref{theorem_semi}), or exponential (Theorem \ref{theorem_exponential}) convergence rates. Our convergence bounds are straightforward to compute, as all elements (e.g., pre-multipliers and exponential rates) are explicitly defined in terms of simple constants and one-step transition expectations.
    \item We devise novel techniques to establish CDs in various scenarios, effectively capturing special dynamics (e.g., reflected boundaries in queueing systems) and parametric dependencies (e.g., traffic intensity in queueing systems and step-size in stochastic algorithms).
    \item We apply our results to constant step-size stochastic gradient descent (SGD) under non-standard assumptions, including infinite variance gradient noise and non-strongly-convex objectives. Our bounds are parametric in the degree of heavy-tailedness of the gradient noise and the degree of flatness of the objective around the optimizer. This analysis sheds light on how these features affect convergence rates (see Section \ref{section_sgd}).
    \item We also apply our results to the G/G/1 queue in heavy traffic as well as stochastic fluid networks. For the G/G/1 queue, we derive a sharp polynomial convergence bound that is uniform with respect to the heavy traffic parameter (see Section \ref{section_lmt}).
    For tandem stochastic fluid networks and related systems, we derive sharp exponential convergence bounds with insightful pre-multipliers (see Section \ref{section_brt}).
    \item One innovative aspect of our analysis is the use of induced metrics. These metrics can be visualized as the minimization of a certain action integral over paths that connect any two given points in the metric space. Thanks to induced metrics, our convergence bounds take an explicit form, involving, for example, a suitable induced metric between the first step of the chain and its initial location. Moreover, these metrics can be used to overcome expansiveness (see Section \ref{section_expansive}).
\end{enumerate}

\subsection{Related works}
To place our contributions in context, we now review the related literature. The findings in \cite{steinsaltz1999locally} bear the closest resemblance to our results. In \cite{steinsaltz1999locally}, a modified transition kernel is introduced to describe a notion of local contraction, which relaxes the global contraction property. This approach yields straightforward convergence bounds but can only be applied to geometrically convergent chains in Euclidean space. In contrast, our CD framework can be applied to sub-geometrically convergent chains in general metric spaces. For geometrically convergent chains, our bound (Theorem \ref{theorem_exponential}), which leverages induced metrics, is stronger than the one in \cite{steinsaltz1999locally}. This is further discussed after Theorem \ref{theorem_exponential}. In \cite{steinsaltz1999locally}, Markov chains are mainly viewed as iterative function systems (IFS), which turns out to be beneficial. For a recent comprehensive survey on IFS, see \cite{ghosh2022iterated}.

In \cite{qin2022geometric}, a bivariate version of \cite{steinsaltz1999locally} is introduced. While \cite{steinsaltz1999locally} enforces the control of drift and contraction point by point, \cite{qin2022geometric} enforces it pair (of points) by pair.
Essentially, in the above two papers, a Lyapunov function (that creates drift) is introduced to modify the original metric, making the Markov chain globally contractive under the modified metric. In \cite{eberle2019quantitative}, it is a concave functional of the original metric that is modified by a Lyapunov function to establish global contraction for geometrically convergent chains.
In the CD framework, we can have two functions: one for metric modification and the other for drift construction.
The two functions are linked via a single inequality (CD).

There is another way to address drift and contraction separately.
The drift and contraction (D\&C) method, starting from \cite{hairer2011asymptotic}, combines contraction inside a selected region and drift outside that region to establish geometric convergence bounds (see, e.g., \cite{jarner2001locally,durmus2015quantitative,douc2018markov}). A recent representative example of this method is Corollary 2.1 in \cite{qin2022geometric}, which also has a random-mapping-representation version in \cite{qin2022wasserstein}. While the D\&C method is capable of handling non-globally-contractive chains, it does face some limitations, as it still requires strict contraction within a selected region. When one uses the D\&C method, there is an implicit trade-off that is captured by the size of the region. Typically, effective contraction requires the region to be small, while effective drift requires the region to be large. Consequently, a suitable choice of region that can generate a sharp bound may not exist (see the discussion at the end of Section \ref{section_sgd}). Balancing drift and contraction becomes even more difficult if one wishes to find bounds that are parametrically sharp across regimes of interest, such as a sequence of queues in heavy traffic. The CD framework turns out to be accurate enough to develop such bounds, as we illustrate in Section \ref{section_lmt} where we derive a polynomial convergence bound that is uniform in heavy traffic.

The D\&C method for polynomially convergent chains is studied in \cite{butkovsky2014subgeometric} and \cite{durmus2016subgeometric}, where the drift conditions are similar to those in \cite{jarner2002polynomial} and \cite{douc2004practical} for estimating polynomial convergence in TV distance. 
They assume that the metric is bounded and the chain is non-expansive.
Their bounds are qualitative in nature, as they are not explicit or might be difficult to compute explicitly. This is further discussed after Theorem \ref{theorem_polynomial_b}.
In the CD framework, quantitative bounds are derived without those assumptions.

Many of the aforementioned results require the chain to be non-expansive, which may not be satisfied in practice. However, it is possible to make the chain non-expansive by modifying the underlying metric. This metric modification approach has been systematically developed for diffusion processes (\cite{eberle2011reflection, eberle2016reflection, zimmer2017explicit, eberle2019quantitativeH, eberle2019sticky}). As we demonstrate in Section \ref{section_expansive}, metric modification and drift construction are naturally integrated under the CD framework.

In summary, the D\&C method modifies the metric to enforce contraction and finds a Lyapunov function to create drift, while the method in \cite{steinsaltz1999locally} and \cite{qin2022geometric} smoothly combines the two steps via a single ``weight'' function. In the current paper, we develop this idea of smooth combination into a unified convergence analysis framework.

The rest of the paper is organized as follows: In Section \ref{section_preliminaries}, we introduce various concepts, including induced metrics and the local Lipschitz constant.
In Section \ref{section_main}, we introduce the contractive drift condition (CD) and present our primary findings, namely, Wasserstein convergence theorems with various convergence rates. In Section \ref{section_tvw}, we highlight several advantages of analyzing convergence in Wasserstein distance over TV distance. In Sections \ref{section_expansive} and \ref{section_sgd}, we use the CD framework to bound the convergence of stochastic algorithms, which can be non-contractive or even expansive. In Sections \ref{section_lmt} and \ref{section_brt}, we introduce two techniques to establish CDs, and we use them to bound the convergence of the G/G/1 queue in heavy traffic, and also stochastic fluid networks. In Section \ref{section_deep}, we describe how the CD framework can allow convergence analysis to be combined with deep learning. All proofs are in Section \ref{section_proofs}.

\section{Preliminaries}
\label{section_preliminaries}
Let $(\X,d)$ be a complete metric space. A curve in $\X$ is a continuous function $\gamma:[0,1]\rightarrow\X$. Given $t\in[0,1]$, the length of $\gamma|_{[0,t]}$ (the restriction of $\gamma$ to $[0,t]$) is given by
$$L(\gamma|_{[0,t]})\defeq\sup_{0=t_0<t_1...<t_n=t,\;n\geq1}\;\;\sum_{k=1}^nd(\gamma(t_{k-1}),\gamma(t_k)).$$
A curve $\gamma$ is rectifiable if $L(\gamma)\defeq L(\gamma|_{[0,1]})<\infty$.
The following path connectivity assumption will be maintained throughout the remainder of this paper.
\begin{assumption}
\label{assumption_path_connect}
Each pair of points in $\X$ is connected by a rectifiable curve.
\end{assumption}
Given a rectifiable curve $\gamma$, its length function $L(\gamma|_{[0,t]})$ is continuously increasing (see, e.g., Chapter 2.3.2 of \cite{burago2001course}), so it induces a finite Borel measure on $[0,1]$.
Given a Borel-measurable function $g:\X\rightarrow\R_+$, the line integral of $g$ along $\gamma$ is well-defined as a Lebesgue-Stieltjes integral (see, e.g., Chapter 6.3.3 of \cite{stein2009real}), namely
$$L(\gamma;g)\defeq\int_0^1 g(\gamma(t))dL(\gamma|_{[0,t]}).$$
If $g$ is bounded away from zero ($\inf_{x\in\X}g(x)>0$), then $g$ induces a metric
$$d_g(x,y)\defeq\inf_{\gamma\in\Gamma(x,y)} L(\gamma;g),\;\;x,y\in\X$$
where $\Gamma(x,y)$ is the set of all rectifiable curves joining $x$ and $y$.
If $g\equiv1$, then $d_g$ is known as the intrinsic metric $d_I$. Under Assumption \ref{assumption_path_connect}, $(\X,d_I)$ is complete; see \cite{hu1978local}. 
In Euclidean space, if $\X$ is a convex set, then $d_I=d$.

Let $X=(X_n:n\geq0)$ be a Markov chain on $\X$ with random mapping representation
$$X_{n+1}=f_{n+1}(X_n),\;\;n=0,1,2,\dots$$
where $f_{n+1}$'s are iid copies of $f$, a random mapping from $\X$ to itself. (In this paper, $n$ is always integer-valued.) In general, a given Markov chain can have many random mapping representations, so our convergence bounds depend upon the particular representation chosen. Starting from initial distribution $X_0$, let
$$X_n\defeq(f_n\circ...\circ f_1)(X_0)\;\;\;\;\text{and}\;\;\;\;\bar{X}_{n}\defeq(f_1\circ...\circ f_n)(X_0)$$
be the forward chain and the backward chain, respectively. For each $n$, $X_n$ and $\bar{X}_n$ have the same marginal distribution, as do their limits if they exist. When the stationary distribution exists, let $X_\infty$ be a random variable having that distribution.
\begin{remark}
\label{remark_notation}
A more commonly used notation for the random mapping representation is $X_{n+1}=f_{\theta_{n+1}}(X_n)$ where $\{f_\theta:\theta\in\Theta\}$ is a functional family and $\theta_{n+1}$'s are iid random variables (\cite{diaconis1999iterated}). In this paper, we write $f_{n+1}(x)$, a univariate random function, instead of $f_{\theta_{n+1}}(x)$, a bivariate deterministic function with a random parameter, because not only $f_{n+1}(x)$ is notationally simpler than $f_{\theta_{n+1}}(x)$ but also $Df_{n+1}(x)$ is simpler than $D_{x}f_{\theta_{n+1}}(x)$ when ``differentiating''.
\end{remark}

The local Lipschitz constant of $f$ at $x\in\X$ is defined as
$$Df(x)\defeq\lim_{\delta\rightarrow0}\sup_{x',x''\in B_\delta(x),\;x'\neq x''}\frac{d(f(x'),f(x''))}{d(x',x'')}$$
where $B_\delta(x)\defeq\{x'\in \X:d(x',x)<\delta\}$. 
In Euclidean space, if $f$ is differentiable, then
$Df(x)=\norm{\nabla f(x)}$, which is the spectral norm of the Jacobian.
The local Lipschitz constant locally describes how expansive or contractive $f$ is around $x$. The following local Lipschitzness assumption will be maintained throughout the remainder of this paper.
\begin{assumption}
\label{assumption_local_lipschitz}
With probability $1$, $f$ is locally Lipschitz, i.e., $Df<\infty$ everywhere.
\end{assumption}
Note that we only assume $Df<\infty$ but not $Df<1$, so $f$ can be expansive (see Section \ref{section_expansive}). Next, we recall the definition of the Wasserstein distance.
Let $\mathcal{P}(\X)$ be the set of integrable probability measures on $\X$ equipped with its Borel sigma-algebra. The Wasserstein distance (induced by $d$) between $\mu,\nu\in\mathcal{P}(\X)$ is
$$W_d(\mu,\nu)\defeq\inf_{\pi\in\mathcal{C}(\mu,\nu)}\int_{\X\times\X}d(x,y)\pi(dx,dy)$$
where 
$$\mathcal{C}(\mu,\nu)\defeq\{\pi\in\mathcal{P}(\X\times\X):\pi(\cdot,\X)=\mu(\cdot),\;\pi(\X,\cdot)=\nu(\cdot)\}$$ is the set of all couplings of $\mu$ and $\nu$.
Given two random variables $Y$ and $Z$, we use $W_d(Y,Z)$ to denote the Wasserstein distance between their marginal distributions. To simplify notation, the Wasserstein distance induced by $d_g$ is denoted by $W_g(\cdot,\cdot)$, and the Wasserstein distance induced by $d_I$, the intrinsic metric, is denoted by $W_I(\cdot,\cdot)$. When $g\geq\epsilon>0$, we have $d_g\geq\epsilon d_I\geq \epsilon d$ and hence $W_g(\cdot,\cdot)\geq \epsilon W_I(\cdot,\cdot)\geq\epsilon W_d(\cdot,\cdot)$. In this paper, we mainly develop upper bounds for $W_I(X_n,X_\infty)$ that are also upper bounds for $W_d(X_n,X_\infty)$. In Euclidean space, if $\X$ is a convex set, then $W_I(\cdot,\cdot)=W_d(\cdot,\cdot)$. In this case, we simply write $W(\cdot,\cdot)$.
\begin{remark}
\label{remark_first_induce}
As far as we are aware, the current paper is the first to use the intrinsic metric $d_I$, induced metric $d_g$, and their corresponding Wasserstein distances to quantify the convergence of Markov chains on general metric spaces. 
In \cite{stenflo2012survey}, $d_I$ and $d_g$ are mentioned, but the author does not derive convergence bounds under these metrics.
\end{remark}

\section{Main results}
\label{section_main}
Given a Markov chain $X$ driven by random mapping $f$, the contractive (transition) kernel is defined as
$$Kh(x)\defeq\E Df(x)h(f(x)),\;\;x\in\X,\;\;h:\X\rightarrow\R_+.$$
Compared with the standard transition kernel given by
$$Ph(x)\defeq\E_x h(X_1)=\E[h(X_1)|X_0=x]=\E h(f(x)),$$ the contractive kernel incorporates the local contraction/expansion information quantified by $Df(x)$.
In \cite{steinsaltz1999locally}, $KV\leq rV$ with $r\in(0,1)$ and $V:\X\rightarrow[1,\infty)$ is used to derive simple geometric convergence bounds for Markov chains in Euclidean space. Now we introduce the contractive drift condition (CD) in general.
Let $\V$ be the family of Borel-measurable functions on $\X$ with positive infima.
\begin{condition}
For $U,V\in\V$, the contractive drift condition $KV\leq V-U$ is
$$KV(x)=\E Df(x)V(f(x))\leq V(x)-U(x),\;\;x\in\X.$$
\end{condition}
In the D\&M or D\&C method, the traditional drift condition ($PV\leq V-U$) can not hold everywhere, because the Lyapunov function can not be indefinitely reduced by the chain. In the region where the drift stops, another condition (minorization or contraction) must be introduced. In contrast, our contractive drift condition can hold everywhere, thanks to the local Lipschitz constant in $K$, which makes it possible to establish convergence in one step.
Now we present the polynomial convergence theorem, the main result of this paper.
The main result is followed by a proof sketch to highlight its novelty.
All proofs are in Section \ref{section_proofs}. Proofs for the current section are in Section \ref{proof_section_main}. Recall that $d_V$ is the metric induced by $V$ while $W_U$ is the Wasserstein distance induced by $d_U$.
\begin{theorem}
\label{theorem_polynomial_b}
Assume that $KU\leq U$ and $KV\leq V-U^{1/b}V^{1-1/b}$ where $U,V\in\V$ and $b>1$. If $\E d_V(X_0,X_1)<\infty$, then $X$ has a unique stationary distribution $X_\infty$ with 
$$W_U(X_n,X_\infty)\leq\sbk{\prod_{k=1}^{\lceil b\rceil-1}\frac{b}{n+k}\cdot\frac{\lceil b\rceil-k}{b-k}}^{\frac{b-1}{\lceil b\rceil-1}}\cdot\E d_V(X_0,X_1).$$
Moreover, $W_U(X_n,X_\infty)=o(1/n^{b-1})$.
\end{theorem}
\begin{proof}[Proof sketch of Theorem \ref{theorem_polynomial_b}]
There are four main steps.
\begin{enumerate}[label=\roman*)]
\item From $KV\leq V-U^{1/b}V^{1-1/b}$, we explicitly extract a CD sequence
$KV_k\leq V_k-V_{k+1}$ where $k=0,\dots,\lceil b\rceil-1$. All $V_k$'s have simple expressions. 
In \cite{jarner2002polynomial,douc2004practical,butkovsky2014subgeometric,durmus2016subgeometric}, a sequence of drift conditions is also extracted from a special drift condition to establish subgeometric convergence bounds, but their extraction cannot be done explicitly, mainly because
the traditional drift condition cannot hold everywhere.
\item Given the explicit CD sequence and some combinatorial identity, we run the chain forward to establish
\begin{equation}
\label{equation_step_two}
    \sum_{n=0}^{\infty}c_n\E DF_n(x)U(F_n(x))\leq V(x),\;\;x\in\X
\end{equation}
where $F_n\defeq f_n\circ\dots\circ f_1$ and $c_n$'s are explicit constants. This step illustrates that the sub-geometric case is harder than the geometric case $KV\leq rV$ (e.g., \cite{steinsaltz1999locally}), where we simply have $\E DF_n(x)V(F_n(x))\leq K^n V(x)\leq r^n V(x).$
\item 
Now we run the chain backward, replacing each $F_n$ in \eqref{equation_step_two} with $\bar{F}_n\defeq f_1\circ\dots\circ f_n$.
Considering all rectifiable curves joining $X_0$ and $f(X_0)$, we integrate \eqref{equation_step_two} and minimize the integral to obtain
$$\sum_{n=0}^\infty c_n\E d_U(\bar{F}_n(X_0),\bar{F}_{n+1}(X_0))=\sum_{n=0}^\infty c_n\E d_U(\bar{F}_n(X_0),\bar{F}_{n}(f(X_0)))\leq \E d_V(X_0,f(X_0))$$
where induced metrics $d_U$ and $d_V$ play a crucial role. On the RHS of \eqref{equation_step_two}, we may integrate $V(x)$ along some simple curve joining $X_0$ and $f(X_0)$. However, on the LHS of \eqref{equation_step_two}, the integral of $D\bar{F}_n(x)U(\bar{F}_n(x))$ along the same curve corresponds to an integral of $U(x)$ along a potentially complicated curve joining $\bar{F}_n(X_0)$ and $\bar{F}_{n+1}(X_0)$. By minimizing the line integral between them, we arrive at the above clean inequality with $d_U$ on the LHS and $d_V$ on the RHS.
\item Combining this inequality with the completeness of $(\X,d_I)$, the backward chain converges a.s. to some $\bar{X}_\infty$. Finally,
$$c_nW_U(X_n,X_\infty)\leq c_n\sum_{k=n}^\infty\E d_U(\bar{X}_k,\bar{X}_{k+1})\leq \sum_{k=n}^\infty c_k\E d_U(\bar{X}_k,\bar{X}_{k+1})\leq\E d_V(X_0,X_1).$$
\end{enumerate}
\end{proof}
A key aspect that distinguishes this bound from the rest of the literature is that it involves an explicit constant and a one-step transition expectation. We believe that this makes it convenient for practical use. This final expression follows as a consequence of our techniques. In contrast, as mentioned in the proof sketch, it is difficult (or at least not direct) to use either D\&M or D\&C methods to obtain a polynomial bound that is computable only in terms of one-step transition expectations.
For example, the main result (Theorem 3) in \cite{andrieu2015quantitative}, while fully characterized, is either difficult to translate into an explicit bound or intended to be qualitative in nature. As a consequence, the corollaries that they provide to derive bounds from Theorem 3 only state the existence of constants. Similarly, the Wasserstein bounds in \cite{durmus2016subgeometric} are intended to be qualitative in nature (the various constants are not explicitly given); see for example the statement of the main result, also labeled Theorem 3. 
In contrast, the only three elements in our assumption ($U,V$, and $b$) directly correspond to the only three terms in our bound ($W_U(X_n,X_\infty)$, $\E d_V(X_0,X_1)$, and the $b$-constant). The constants that we provide can be further simplified to facilitate their evaluation. For instance, $W_U(X_n,X_\infty)$ is lower bounded by $\inf_{x\in\X}U(x)\cdot W_I(X_n,X_\infty)$ while $d_V(X_0,X_1)$ is upper bounded by any $V$-integral from $X_0$ to $X_1$.

Since the polynomial bound in Theorem \ref{theorem_polynomial_b} is as simple as the exponential bound in Theorem \ref{theorem_exponential} (below), one may wonder whether the assumption in Theorem \ref{theorem_polynomial_b} is as strong as the assumption in Theorem \ref{theorem_exponential} (e.g., global contraction/non-expansion). Fortunately, the answer is no. At first glance, $KU\leq U$, which implies $\E d_U(f(y),f(z))\leq d_U(y,z)$, may look like a non-expansion assumption. However, we have the freedom to carefully design a non-trivial $U$ that makes an expansive chain non-expansive under the $U$-induced metric, which aligns with the theme of metric modification in the recent literature on Wasserstein convergence (e.g., \cite{eberle2019quantitative}). The application of the CD framework to an expansive chain is in Section \ref{section_expansive}.
If the chain is already non-expansive under the original metric (e.g., G/G/1), then $KU\leq U$ holds with $U\equiv1$, so $W_U=W_1$ becomes $W_I$ the Wasserstein distance induced by the intrinsic metric $d_I$.
The following corollary is convenient in practice.
\begin{corollary}
\label{corollary_polynomial_m}
Assume that $K\mathbf{1}\leq\mathbf{1}$ and $KV\leq V-\delta V^{1-1/m}$ where $V\in\V$, $\delta>0$, and integer $m\geq2$. If $\E d_V(X_0,X_1)<\infty$, then $X$ has a unique stationary distribution $X_\infty$ with 
$$W_I(X_n,X_\infty)\leq\frac{1}{\delta^m}\sbk{\prod_{k=1}^{m-1}\frac{m}{n+k}}\cdot\E d_V(X_0,X_1).$$
\end{corollary}
The explicit polynomial bound in Theorem \ref{theorem_polynomial_b} also allows us to optimize $b$ for each $n$ when there is a range of available $b$'s, which leads to bounds with other rates. We illustrate this by establishing an explicit semi-exponential bound (e.g., $\exp(-\sqrt{n})$). Its D\&M and D\&C counterparts can be found in \cite{douc2004practical} and \cite{durmus2016subgeometric}, respectively. 
\begin{theorem}
\label{theorem_semi}
Assume that $K\mathbf{1}\leq\mathbf{1}$ and $KV\leq V-\delta V/(\log V)^\lambda$ where $\delta,\lambda>0$, $V\in\V$, and $V>1$. If $\E d_V(X_0,X_1)<\infty$, then $X$ has a unique stationary distribution $X_\infty$ with
$$W_I(X_n,X_\infty)\leq e^{1/\eta}(n/2)\exp({-n^\eta\kappa^\eta/(e\eta)})\cdot\E d_V(X_0,X_1),\;\;n\geq (2e)^{1/\eta}/\kappa$$
where $\eta=1/(\lambda+1)$ and $\kappa=\delta(e/\lambda)^\lambda$.
\end{theorem}
\begin{proof}[Proof sketch of Theorem \ref{theorem_semi}]
For integer $m\geq2$, we have $$KV\leq V-\delta V/(\log V)^\lambda\leq V-(\kappa/m^\lambda)V^{1-1/m}$$
and the corresponding polynomial bound given by Corollary \ref{corollary_polynomial_m}, which becomes the above semi-exponential bound when $m=cn^{\eta}$ with some $c>0$.
\end{proof}
Last, but not least, we present the geometric bound in the CD framework.
\begin{theorem}
\label{theorem_exponential}
Assume that $KV\leq rV$ where $V\in\V$ and $r\in(0,1)$. If $\E d_V(X_0,X_1)<\infty$, then $X$ has a unique stationary distribution $X_\infty$ with 
$$W_V(X_n,X_\infty)\leq[r^n/(1-r)]\cdot\E d_V(X_0,X_1).$$
\end{theorem}
\begin{proof}[Proof sketch of Theorem \ref{theorem_exponential}]
As mentioned in the previous proof sketch, the geometric case is much simpler than the sub-geometric case.
Considering all rectifiable curves joining $X_0$ and $f(X_0)$, we integrate $\E D\bar{F}_n(x)V(\bar{F}_n(x))\leq K^n V(x)\leq r^n V(x)$ and minimize the integral to obtain
\begin{align*}
    &W_V(X_n,X_\infty)\\
    \leq& \sum_{k=n}^\infty\E d_V(\bar{F}_k(X_0),\bar{F}_{k+1}(X_0))=\sum_{k=n}^\infty\E d_V(\bar{F}_k(X_0),\bar{F}_k(f(X_0)))\leq \sum_{k=n}^\infty r^k\cdot\E d_V(X_0,X_1).
\end{align*}
\end{proof}
In Euclidean space, a similar bound can found in \cite{steinsaltz1999locally}, namely
$$W(X_n,X_\infty)\leq[r^n/(1-r)]\cdot\E \sbk{\norm{X_1-X_0}\cdot\sup_{t\in[0,1]}V((1-t)X_0+tX_1)}$$
where $V\geq1$, $\norm{\cdot}$ is the Euclidean norm, and $W(\cdot,\cdot)$ is the corresponding Wasserstein distance. This bound is weaker than Theorem \ref{theorem_exponential} because $W(X_n,X_\infty)$ on the LHS is upper bounded by $W_V(X_n,X_\infty)$ while the expectation on the RHS is lower bounded by $\E d_V(X_0,X_1)$. When the chain is globally contractive, i.e., there exists $r\in(0,1)$ such that $\E d(f(y),f(z))\leq r d(y,z)$ for all $y,z\in\X$, it is well known that
$W(X_n,X_\infty)\leq[r^n/(1-r)]\cdot\E d(X_0,X_1)$ (see, e.g., \cite{stenflo2001markov}).
In the CD framework, the global contraction means $\E Df(x)\leq r$ for all $x\in\X$. Applying Theorem \ref{theorem_exponential} with $V\equiv1$ yields
$W_I(X_n,X_\infty)\leq[r^n/(1-r)]\cdot\E d_I(X_0,X_1).$ The two bounds have the same form but the original metric $d$ in the former is replaced by the intrinsic metric $d_I$ in the latter. This is because integrating $\E Df(x)\leq r$ along some path leads to $\E d_I(f(y),f(z))\leq r d_I(y,z)$ but not $\E d(f(y),f(z))\leq r d(y,z)$. In Euclidean space, if $\X$ is a convex set, then we do not need to distinguish between $d$ and $d_I$.

\section{From total variation to Wasserstein}
\label{section_tvw}
Before diving into examples, we briefly compare the TV distance and the Wasserstein distance. In the literature, the TV distance has long been the standard metric used to measure the convergence of Markov chains. Here, we use a popular example to illustrate that the Wasserstein distance can be a better choice.

Nowadays, machine learning models are trained on large but finite datasets. To minimize the loss over the whole dataset, stochastic gradient descent (SGD) is widely used because computing the exact gradient is too expensive. When the step-size is constant, SGD is a time-homogeneous Markov chain. Since SGD samples from a finite dataset, the support of its transition kernel is discrete. Therefore, it is unrealistic to assume that the transition kernel has a continuous component, which is required to establish a minorization condition (e.g., \cite{yu2021analysis}). In fact, SGD typically does not converge in TV distance at all. For example, to solve
$$\min_{x\in\R^d}\E\norm{x-Y}^2/2$$
where $Y$ is a discrete random vector and $\norm{\cdot}$ is the Euclidean norm, the SGD iteration with step-size $\alpha\in(0,1)$ is
$$X_{n+1}=X_n-\alpha(X_n-Y_{n+1})=(1-\alpha)X_n+\alpha Y_{n+1}$$
where $Y_{n+1}$'s are iid copies of $Y$. This is just an AR(1) process. When $Y$ is not constant, it is known that the stationary distribution of an AR(1) process can not contain a point mass; see, e.g., Proposition 2.5.2 in \cite{buraczewski2016stochastic}. Starting from a fixed point, the $n$-step marginal distribution is always discrete, so it cannot converge to its atomless limit in TV distance. However, it converges in Wasserstein distance because of the global contraction. Moreover, the contraction rate $1-\alpha$ (for Wasserstein convergence) is dimension-free while there is no minorization condition (for TV convergence). 
In high-dimensional spaces, even if a minorization condition can be established (e.g., when $Y$ is normally distributed), the resulting TV convergence rate scales poorly with dimension (\cite{qin2022wasserstein}).

\section{Application to expansive ULA}
\label{section_expansive}
The contractive drift framework can handle locally expansive or even non-locally expansive chains. The name emphasizes the collaborative contribution of contraction ($Df<1$) and drift $(PV<V)$ to the convergence, but expansion ($Df>1$) and anti-drift ($PV>V$) are allowed , as justified by the following example.

Suppose that we want to sample from the following distribution that has a semi-exponential tail
\begin{equation*}
\pi(x)=e^{-g(x)},\quad
g(x)=\left\{
        \begin{array}{ll}
        \sqrt{\abs{x}} & \quad |x|\geq L \\
        ax^2+c & \quad |x|<L
        \end{array}
    \right.
\end{equation*}
where $L>0$ and $a,c>0$ make $g,g'$ continuous at $\pm L$. In particular, $a=1/(4L^{3/2})$.
The random mapping representation of the corresponding unadjusted Langevin algorithm (ULA) with step-size $\gamma>0$ is
\begin{equation}
\label{def_ula}
f(x) = x-\gamma g'(x)+\sqrt{2\gamma}Z=\left\{
        \begin{array}{ll}
        x-\frac{\text{sgn}(x)\gamma}{2\sqrt{|x|}}+\sqrt{2\gamma}Z & \quad |x|\geq L \\
        x-2\gamma ax+\sqrt{2\gamma}Z & \quad |x|<L
        \end{array}
    \right.
\end{equation}
where $\text{sgn}(x)=I(x>0)-I(x<0)$ and $Z\sim N(0,1)$. The local Lipschitz constant is
\begin{equation*}
Df(x) =\left\{
        \begin{array}{ll}
        1+\frac{\gamma}{4|x|^{3/2}} & \quad |x|\geq L\\
        1-2\gamma a & \quad |x|<L
        \end{array}
    \right.,
\end{equation*}
which shows that the Markov chain is not contractive. It is {\it non-locally} expansive as $Df(x)>1$ when $|x|\geq L$. 
In the literature of ULA, most works focus on the geometrically convergent case (e.g., \cite{durmus2017nonasymptotic,durmus2019high}).
However, we should not expect the above ULA to converge geometrically fast, as its drift $-\gamma/(2\sqrt{x})$ vanishes as $x\gti$, making it hard to establish a geometric drift condition, let alone the expansion.
Before presenting the polynomial convergence bound for this Markov chain, we introduce some notations: $x\wedge y=\min(x,y)$, $x\vee y=\max(x,y)$, and $y_n=\Theta(x_n)$ means $0<\liminf_{n\gti}(y_n/x_n)<\limsup_{n\gti}(y_n/x_n)<\infty$.
Proofs for the current section are in Section \ref{proof_section_expansive}.
\begin{proposition}
\label{proposition_expansive}
Let $X$ be the Markov chain defined by \eqref{def_ula}. If $L^{1/4}$ is an even integer larger than $4$ and $\gamma\leq 2^{13-2L^{1/4}}/(\E Z^{L^{1/4}})^2$, then $X$ has a unique stationary distribution $X_\infty$ with
\begin{align*}
W(X_n,X_\infty)\leq\frac{1}{(\gamma/56)^b4L^{3/2}}\cdot\sbk{\prod_{k=1}^{\lceil b\rceil-1}\frac{b}{n+k}\cdot\frac{\lceil b\rceil-k}{b-k}}^{\frac{b-1}{\lceil b\rceil-1}}\cdot\E\sbk{\int_{X_0\wedge X_1}^{X_0\vee X_1}V(x)dx}
\end{align*}
where $b=(2/3)(L^{1/4}-2)$ and $V(x)=x^{L^{1/4}}+(5/2)L^{L^{1/4}}$. In particular, $W(X_n,X_\infty)=O(1/n^{b-1})$ as $n\gti$ where $b-1=\Theta(L^{1/4})$ as $L\gti.$
\end{proposition}
To apply Theorem \ref{theorem_polynomial_b}, the first step is to show that $KU\leq U$ (metric modification) where $U(x)=x^2+4L^{3/2}$. The second step is to establish that $KV\leq V-(\gamma/56)U^{1/b}V^{1-1/b}$ (drift construction) where $V(x)=x^m+M$. The two parameters $m$ and $M$ are carefully chosen to balance the drift in the expansive region and the anti-drift in the contractive region, which leads to the parametric polynomial convergence rate $1/n^{(2/3)(L^{1/4}-7/2)}$. Here $L^{1/4}$ is assumed to be an even integer to simplify the analysis of $V$.
The expansive ULA example illustrates how metric modification and drift construction are naturally integrated under the CD framework to generate quantitative bounds. For more complex chains where contraction does not hold ``in the bulk'' (e.g., ULA sampling from multimodal distributions), establishing CDs is theoretically possible, but
handpicking a suitable $U$ to fully eliminate expansion between contractive regions is practically challenging. Deep learning may be leveraged to construct $U$; see Section \ref{section_deep}.

The polynomial bounds in \cite{butkovsky2014subgeometric} and \cite{durmus2016subgeometric}, despite being qualitative in nature, are not applicable here because they require the metric to be bounded by one and do not allow expansion.
Although a new metric (e.g., the $U$-induced metric above) can be constructed to eliminate expansion, it is not natural to enforce a bounded metric when the algorithm explores the whole Euclidean space.

\section{Application to non-standard SGD}
\label{section_sgd}
Stochastic gradient descent (SGD) and its variants have achieved remarkable empirical success in training neural networks, which may be attributed to their ability to find flat minima in the loss landscape that lead to better generalization. During the training process, heavy-tailed gradient noise is often observed, and it is used to explain why SGD tends to prefer flat minima (\cite{simsekli2019tail}). Under the CD framework, we explicitly bound the convergence of stylized {\it non-standard} SGD to understand how its performance is affected by the degree of flatness of the minima (e.g., quartic basin) and the degree of heavy-tailedness of the gradient noise (e.g., infinite variance). 
The global contraction result under the standard assumption (strongly convex objective, Lipschitz gradient, finite-variance gradient noise) can be found in \cite{dieuleveut2020bridging}.
Powered by deep learning, the CD framework can also be applied to realistic SGD; see \cite{qu2024deep} for an example.
Proofs for the current section are in Section \ref{proof_section_sgd}.

\subsection{Non-strongly-convex objective}
Suppose that we want to minimize the following non-strongly-convex objective that is flat around the origin
\begin{equation*}
h(x) = \left\{
        \begin{array}{ll}
        |x|^m/m & \quad |x| < 1 \\
        x^2/2-1/2+1/m & \quad |x|\geq1
        \end{array}
    \right.,
\end{equation*}
where $m\geq3$.
The stylized SGD iteration with step-size $\alpha\in(0,1)$ and iid unbiased gradient noise $Z$ is
\begin{equation}
\label{def_sgd_1}
    f(x)=x-\alpha(h'(x)+Z)= \left\{
        \begin{array}{ll}
        x-\alpha(\text{sgn}(x)|x|^{m-1}+Z) & \quad |x| < 1 \\
        x-\alpha(x+Z) & \quad |x|\geq1
        \end{array}
    \right..
\end{equation}
Compared with ULA, here $Z$ may not be normally distributed, and it is multiplied by $\alpha$ rather than $\sqrt{2\alpha}$.
The local Lipschitz constant is
\begin{equation*}
Df(x) =\left\{
        \begin{array}{ll}
        1-\alpha(m-1)|x|^{m-2} & \quad |x|<1\\
        1-\alpha & \quad |x|\geq1
        \end{array}
    \right..
\end{equation*}
Since $Df(0)=1$, the chain is not globally contractive. We use a wedge-like function $V(x)=\delta(1-|x|)_++1$ to artificially create drift/contraction around the origin. Since $V$ reaches it maximum at the origin, we must have $KV(0)<V(0)$.

\begin{proposition}
\label{proposition_sgd_1}
    Let $X$ be the Markov chain defined by \eqref{def_sgd_1}.
    If $\alpha\leq(1/8)/\E(1+|Z|)$, then $X$ has a unique stationary distribution $X_\infty$ with
    $$W(X_n,X_\infty)\leq(2/\tilde{\alpha}^{m-2})\cdot\prs{1-\tilde{\alpha}^{m-2}\alpha}^n\cdot \E\abs{h'(X_0)+Z_1}$$
    where $\tilde{\alpha}=(1-\E(1-\alpha|Z|)_+)/4=\Theta(\alpha)$ as $\alpha\gtz.$
    In particular, $W(X_n,X_\infty)=O(r^n)$ as $n\gti$ where $1-r=\Theta(\alpha^{m-1})$ as $\alpha\gtz.$
\end{proposition}
This bound explicitly describes how the flatness near the optimizer affects the convergence when the step-size is small. The larger the power $m$, the flatter the objective $h$, the smaller the gap $1-r$, the slower the convergence.
Note that when $m=2$ the objective becomes $h(x)=x^2/2$. The SGD becomes $f(x)=(1-\alpha)x-\alpha Z$ where the gap is clearly $\alpha$, indicating that $1-r=\Theta(\alpha^{m-1})$ may be sharp. This SGD example illustrates how the metric is modified under the CD framework to establish global contraction. In particular, $V$ with a ``wedge'' at the non-contractive region can restore contraction, because the corresponding $d_V$ locally ``stretches'' the original metric.
For stochastic algorithms exploring complex landscapes, this ``Riemannian'' metric modification (induced metrics) appears to be a suitable approach for restoring contraction, as it addresses local non-contraction in a localized manner. In contrast, \cite{eberle2019quantitative} modify the metric globally by applying a concave function (i.e., $f_0(d(x,y))$ with $f_0$ concave).

\subsection{Heavy-tailed gradient noise}
Suppose that we want to minimize the following objective that is a generalization of the Huber loss
\begin{equation*}
h(x) = \left\{
        \begin{array}{ll}
        x^2/2 & \quad |x| < 1 \\
        |x|^\beta/\beta-1/\beta+1/2 & \quad |x|\geq1
        \end{array}
    \right.,
\end{equation*}
where $\beta\in(1,2).$ 
The stylized SGD iteration with step-size $\alpha\in(0,1)$ and iid unbiased gradient noise $Z$ is
\begin{equation}
\label{def_sgd_2}
    f(x)=x-\alpha(h'(x)+Z)= \left\{
        \begin{array}{ll}
        x-\alpha(x+Z) & \quad |x| < 1 \\
        x-\alpha(\text{sgn}(x)|x|^{\beta-1}+Z) & \quad |x|\geq1
        \end{array}
    \right..
\end{equation}
We assume that $Z$ has a finite $\gamma$-th moment with $\gamma\in(1,2)$. The local Lipschitz constant is
\begin{equation*}
Df(x) =\left\{
        \begin{array}{ll}
        1-\alpha & \quad |x|<1\\
        1-\alpha(\beta-1)|x|^{\beta-2} & \quad |x|\geq1
        \end{array}
    \right.,
\end{equation*}
Since $Df(x)\rightarrow1$ as $x\gti$, the chain is not globally contractive. We should not expect this chain to converge geometrically fast, as its drift $-\alpha|x|^{\beta-1}$ vanishes as $x\gti$, and furthermore $Z$ is heavy-tailed.

\begin{proposition}
\label{proposition_sgd_2}
    Let $X$ be the Markov chain defined by \eqref{def_sgd_2}. If $\beta+\gamma>3$ and $\alpha$ is small enough such that
    $$P(Z\leq1/\alpha-1)\geq3/4,\;\;\sup_{z\geq1/\alpha-1}\E(-Z)I(Z\leq z)\leq1/8,\;\;\alpha^{\gamma-1}\E|Z|^\gamma/(\gamma-1)\leq1/8,$$
    then $X$ has a unique stationary distribution $X_\infty$ with
    \begin{align*}
    W(X_n,X_\infty)\leq\frac{1/\bar{M}^{b-1}}{(\gamma-1)^b(\alpha/2)^b}\cdot\sbk{\prod_{k=1}^{\lceil b\rceil-1}\frac{b}{n+k}\cdot\frac{\lceil b\rceil-k}{b-k}}^{\frac{b-1}{\lceil b\rceil-1}}\cdot\E\sbk{\int_{X_0\wedge X_1}^{X_0\vee X_1}V(x)dx}
    \end{align*}
    where $b=(\gamma-1)/(2-\beta)$, $\bar{M}=\E(1+|Z|)^{\gamma-1}/\alpha+(\gamma+1)/2$, and $V(x)=|x|^{\gamma-1}+\bar{M}-1.$
    In particular, $W(X_n,X_\infty)=O(1/n^{b-1})$ as $n\gti$ where $b-1=(\gamma+\beta-3)/(2-\beta)$.
\end{proposition}
This bound explicitly describes how the heavy-tailedness of the gradient noise ($\gamma$) and the growth rate of the objective ($\beta$) both contribute to the convergence when $\beta+\gamma>3$. When $\beta+\gamma=3$, the SGD may not have a polynomial rate of convergence. For example, when $\beta=1$, the SGD has constant drift toward the origin when it is far away. The waiting time sequence of the G/G/1 queue also has this feature, the convergence of which is studied in Section \ref{section_lmt}. Proposition \ref{proposition_exact} shows that the waiting time sequence 
converges, but not at any polynomial rate, when the noise only has two finite moments ($\gamma=2$).

As mentioned in the Introduction, the implicit trade-off captured by the size of the selected region may prevent the D\&C method from obtaining sharp bounds. In the above SGD example, the larger the selected region, the stronger the drift outside, the weaker the contraction inside. In the D\&C method, the convergence bound is a combination of the worst drift rate outside and the worst contraction rate inside (both of them are reached on the boundary). In contrast, the CD framework allows us to {\it smoothly} combine drift and contraction, so we do not need to compute the two worst rates.

\section{Large M technique and the G/G/1 queue}
\label{section_lmt}
In both Propositions \ref{proposition_expansive} and \ref{proposition_sgd_2}, we consider $V(x)=x^m+M$ and tune $m,M$ to establish polynomial CD. We call this technique the large M technique. In the following, using the waiting time sequence of the G/G/1 queue as an example, we explain the idea behind this technique.

Although we obtain parametric polynomial bounds in Propositions \ref{proposition_expansive} and \ref{proposition_sgd_2}, for these non-standard examples, it is hard to tell whether the parameter dependency is optimal. Given the simplicity of the G/G/1 queue, we are able to rigorously verify that the polynomial bound established under the CD framework is sharp (exact polynomial rate) and parametrically sharp (heavy traffic uniformity). Proofs for the current section are in Section \ref{proof_section_lmt}.

\subsection{Large M technique}
For the waiting time sequence of the G/G/1 queue, the random mapping representation and its local Lipschitz constant are
\begin{equation}
\label{def_gg1}
    f(x)=(x+Z)_+,\;\;Df(x)=I(x+Z\geq0),\;\;x\geq0
\end{equation}
where $Z$, the difference between the service time and the interarrival time, has negative mean. 
When $x+Z<0$, the local Lipschitz constant of $f$ at $x$ is $0$, because $f$ maps not only $x$ but also a small neighborhood around it to a single point, the origin. Let $\delta=-\E Z>0$ and $V(x)=x+1$. We know that $V$ is a traditional Lyapunov function, i.e., $PV\leq V-\delta/2$ for large $x$. Since $Df\leq1$, we immediately have $KV\leq V-\delta/2$ for large $x$. To make the inequality hold everywhere, we can simply add a large constant $M$ to $V$. Now we explain why it works. Suppose that $KV(x)>V(x)-\delta/2$ at some $x$. Then $\E Df(x)=P(x+Z\geq0)<1$ at this $x$, because $P(x+Z\geq0)=1$ implies $PV(x)=V(x)-\delta$. When adding $M$ to $V$,
$$\E Df(x)(V(f(x))+M)-(V(x)+M)=KV(x)-V(x)-(1-\E Df(x))M,$$
which becomes less than $-\delta/2$ when $M$ is large enough. Tuning $M$ is an {\it algebraic} way to balance drift and contraction, which is done geometrically in the D\&C method (region selection). A good choice of $M$ is the key to establish sharp bounds.

\subsection{Exact rate of convergence}
Compared with the total variation distance, a significant feature of the Wasserstein distance is its integrability requirement, i.e., $W_d$ measures the distance between two distributions that are integrable with respect to $d$.
We need to take care in respecting this requirement when modifying $d$.
For example, consider the point mass at the origin ($\delta_0$), random variable $Z$, and function $V(x)=|x|^m$. Since there is only one coupling between them,
$$W_{V}(\delta_0,Z)=\E d_V(0,Z)=\E\int_0^{|Z|}x^m dx=\E\sbk{\frac{x^{m+1}}{m+1}\Big|_0^{|Z|}}=\frac{\E|Z|^{m+1}}{m+1},$$
which is finite if and only if $Z$ has $m+1$ finite moments. In fact, this requirement can guide us in choosing the correct $V$ (e.g., if $\E|Z|^7<\infty$, then $V(x)=|x|^6$). However, this requirement is circumvented in \cite{butkovsky2014subgeometric} and \cite{durmus2016subgeometric} by assuming $d\leq1,$ which explains why their sub-geometric Wasserstein convergence rates are the same as the corresponding TV convergence rates; see Table 1 of \cite{durmus2016subgeometric}. In the following, for the waiting time sequence of the G/G/1 queue, we compute the exact polynomial rate of convergence in the standard Wasserstein distance, which is different from the corresponding TV convergence rate. We use Spitzer's identity (\cite{spitzer1956combinatorial}) to show that the polynomial rate obtained under the CD framework is exact.

\begin{proposition}
\label{proposition_exact}
    Let $X$ be the Markov chain defined by \eqref{def_gg1} starting from $0$. Let $m$ be a positive integer. If $\E Z_+^{m+1}<\infty$ but $\E Z_+^{m+1+\epsilon}=\infty$ for all $\epsilon>0$, then $X$ has a unique stationary distribution $X_\infty$ with
    $$\limsup_{n\gti}n^{m-1}W(X_n,X_\infty)=0,\;\;\;\;\limsup_{n\gti}n^{m-1+\epsilon}W(X_n,X_\infty)=\infty$$
    for all $\epsilon>0.$
\end{proposition}

When $Z_+$ only has $m+1$ finite moments, the limit $X_\infty$ only has $m$ finite moments; see \cite{kiefer1956characteristics}. Since one moment is needed to define the Wasserstein distance, the exact polynomial rate of convergence is naturally $m-1$, which is strictly slower than the corresponding TV convergence rate $m$ (\cite{jarner2002polynomial}). To be specific, $PV\leq V-cV^{1-1/(m+1)}$ with $V(x)=|x|^{m+1}+C$ leads to TV rate $m$, while $KV\leq V-cV^{1-1/m}$ with $V(x)=|x|^m+C$ leads to Wasserstein rate $m-1$, where the power is reduced by one to meet the integrability requirement.

\subsection{Uniform convergence in heavy traffic}
After demonstrating that the CD framework can generate sharp bounds, now we show that it can also generate parametrically sharp bounds. When the G/G/1 queue is in heavy traffic, the random mapping representation of its waiting time sequence becomes 
\begin{equation}
\label{def_ht}
f^\delta(x)=(x+Y-\delta)_+,\;\;x\geq0
\end{equation}
where $Y$ has zero mean and $\delta\downarrow0.$ Let $X^\delta$ be the Markov chain defined by $f^\delta$.
Smaller downward drift $\delta$ means greater congestion in the system, i.e., the system converges slower and is more likely to reach large values. 
However, as long as $\E Y^2<\infty,$ the scaled process $\delta X_{n/\delta^2}^\delta$ can be well approximated by a reflected Brownian motion (\cite{harrison1981reflected}) that converges exponentially fast (\cite{budhiraja2007long}). If a convergence bound for $X^\delta$ is sharp in $\delta$, then the corresponding bound for $\delta X^\delta_{n/\delta^2}$ should not explode as $\delta\downarrow0.$ It turns out that the CD framework can generate a bound with this property.
\begin{proposition}
\label{proposition_uniform}
Let $X^\delta$ be the Markov chain defined by \eqref{def_ht} starting from $0$.
Assume that $\E Y_+^{m+1}<1$ with integer $m\geq1$.
Further assume that $Y$ is not bounded from below and
\begin{equation}
\label{minus_b}
    -b=\inf_{y\in[-1,\infty)}\E\sbk{Y+y\Big|Y+y\leq0}>-\infty.
\end{equation}
Then, $X^\delta$ has a unique stationary distribution $X_\infty^\delta$ with
\begin{align*}
    &\sup_{\delta\in(0,1)}W\prs{\delta X_{n/\delta^2}^\delta,\delta X_\infty^\delta}\\
    \leq&\frac{4}{m}\sbk{\frac{16\E(2+|Y|)^m(1+b)^m}{n}}^{m-1}\cdot\E\sbk{\frac{\prs{1+Y_+}^{m+1}}{m+1}+(1+b)^mY_+}.
\end{align*}
\end{proposition}

Assumption \eqref{minus_b} states that $-Y$ has a bounded residual mean ``lifetime''; that is, conditional on $-Y\geq y$, the expected overshoot $(-Y)-y$ is bounded from above. 
This property is present in exponential distributions but not in Pareto distributions. Although this property does not directly correspond to having light tails, it intuitively makes $-Y$ less likely to take unusually large values.

Although limiting the upper tail of $Y$ (e.g., $\E Y_+^{m+1}<\infty$) is sufficient to establish convergence (Proposition \ref{proposition_exact}), we believe that limiting the lower tail of $Y$ (e.g., \eqref{minus_b}) is necessary to make the convergence uniform in heavy traffic (Proposition \ref{proposition_uniform}). Here, we present an intuitive explanation. Recall that $Df^\delta(x)=I(x+Y-\delta\geq0)$, i.e., contraction ($Df^\delta(x)=0$) only happens at the origin ($f^\delta(x)=0$). When the chain is around the origin, its contraction rate depends on how frequently the origin is visited. 
We compare $\bar{X}^\delta_{n+1}=(\bar{X}^\delta_n+\bar{Y}_{n+1}-\delta)_+$ and $\tilde{X}^\delta_{n+1}=(\tilde{X}^\delta_n+\tilde{Y}_{n+1}-\delta)_+$ where $\bar{Y}\sim N(0,1)$ while $\tilde{Y}$ {\it only} have two finite moments ($\E\tilde{Y}=0$, $\E\tilde{Y}^2=1$). As $\delta\downarrow0$, $\delta\bar{X}^\delta_{n/\delta^2}$ and $\delta\tilde{X}^\delta_{n/\delta^2}$ can be approximated by the same reflected Brownian motion, so they spend a similar proportion of time around the origin, but they may visit the origin at different frequencies. Let $\delta=0.1$ and $\bar{X}_0=\tilde{X}_0=0$. Consider a typical light-tailed sequence and a typical heavy-tailed sequence
$$(\bar{Y}_1,...,\bar{Y}_8)=(3,-3,3,3,-3,-3,3,-3),\;\;(\tilde{Y}_1,...,\tilde{Y}_8)=(1,1,1,1,-7,1,1,1).$$
They have the same sample mean $0$ and similar sample variances. Driven by these two sequences, $\bar{X}^\delta_{n}$ visits the origin three times, while $\tilde{X}^\delta_{n}$ visits the origin only once. This comparison shows that the heavy-tailedness of the lower tail of $Y$ can slow down the contraction. This slow-down effect becomes severer as $\delta\downarrow0$ as time $n$ is accelerated by $(1/\delta^2)$ in the scaled process. Therefore, to establish uniform convergence in heavy traffic, not only the upper tail but also the lower tail of $Y$ should be limited, and \eqref{minus_b} is one way to do so.

\section{Boundary removal technique and stochastic fluid networks}
\label{section_brt}
The large M technique is useful in establishing sub-geometric CDs, but not geometric CDs, because 
$$\E Df(x)(V(f(x))+M)-r(V(x)+M)=KV(x)-rV(x)-(r-\E Df(x))M$$
may not decrease as $M$ increases. Fortunately, for stochastic systems with reflecting boundaries, we have a simple technique to establish geometric CDs, which we call the boundary removal technique. In the following, using one-dimensional reflected Brownian motion (RBM) as an example, we explain the idea behind this technique. 

We use this technique to bound the convergence of tandem stochastic fluid networks and related systems. In the resulting convergence bound, the exponential rate is sharp, and the pre-multiplier provides insight into how the one-step transition structure affects convergence. Proofs for the current section are in Section \ref{proof_section_brt}.

\subsection{Boundary removal technique}
Let $X=(X_t:t\geq0)$ be the RBM solving the following stochastic differential equation (SDE)
\begin{equation}
\label{RBM}
    dX_t=-rdt+\sigma dB_t+dL_t
\end{equation}
where $r,\sigma>0$, $B$ is a standard Brownian motion (BM), and $L$ is a continuous non-decreasing process for which $I(X_t>0)dL_t=0$ for all $t>0$. Starting from $X_0=x\geq0$, by Theorem 6.1 of \cite{chen2001fundamentals}, we have
$$X_t=Z_t+L_t,\;\;Z_t=x-rt+\sigma B_t,\;\;L_t=\sup_{0\leq s\leq t}(-Z_t)_+=\max\prs{0,\sup_{0\leq s\leq t}(-Z_t)}$$
where $Z$ is a ``free'' BM drifting downward while $L$ is the regulator that keeps $X$ non-negative. For $s>0$, let $X^s=(X_{ns}:n\geq0)$ be the $s$-skeleton of $X$. The random mapping representation of $X^s$ is
\begin{equation}
\label{rmr_RBM}
f^s(x)=x-rs+\sigma B_s+\max\prs{0,-x+\sup_{0\leq u\leq s}(-(-ru+\sigma B_u))}.
\end{equation}
The local Lipschitz constant is
$$Df^s(x)=I\prs{\sup_{0\leq u\leq s}(-(-ru+\sigma B_u))\leq x}=I\prs{\inf_{0\leq u\leq s}(x-ru+\sigma B_u)\geq 0}=I(\tau>s)$$
where $\tau=\inf\{t>0:Z_t<0\}$. Similar to the G/G/1 queue, contraction happens only when the origin is visited (by $Z$ during $[0,s]$). Let $K^s$ be the contractive kernel of $X^s$. For any positive function $V$, we have
$$K^sV(x)=\E Df^s(x)V(f^s(x))=\E_x I(\tau>s)V(X_s)=\E_x I(\tau>s)V(Z_s)\leq \E_x V(Z_s)$$
where the third equality holds because the RBM and the free BM are the same until they hit the origin ($s<\tau\;\Rightarrow\;X_s=Z_s$). Now, it suffices to find a drift condition for the free BM, as if the boundary does not exist. Let $V_c(x)=e^{cx}$ with $c>0$. Then
$$\E_x V_c(Z_s)=\E e^{c(x-rs+\sigma B_s)}=V_c(x)e^{-crs+c^2\sigma^2 s/2}=V_c(x)e^{-r^2s/(2\sigma^2)}$$
where $c$ is optimized by $r/\sigma^2.$ By Theorem \ref{theorem_exponential}, we have
$$W(X_{ns},X_\infty)\leq\frac{\lambda^{ns}}{1-\lambda^s}\cdot\E \int_{X_0\wedge X_s}^{X_0\vee X_s}e^{cx}dx=\frac{\lambda^{ns}}{1-\lambda^s}\cdot\frac{\E\abs{e^{cX_s}-e^{cX_0}}}{c}$$
where $\lambda=e^{-r^2/(2\sigma^2)}$. This Wasserstein convergence rate matches the exact TV convergence rate obtained in \cite{glynn2018rate}. 
When $t$ is a multiple of $s$, the above bound becomes $W(X_t,X_\infty)\leq C\lambda^t$.
When $t$ is not a multiple of $s$,
$$W(X_t,X_\infty)=W\prs{f^{t-[t/s]s}(X_{[t/s]s}),f^{t-[t/s]s}(X_\infty)}\leq W(X_{[t/s]s},X_\infty)\leq C\lambda^{t-s}$$
where the first inequality is because $f$ is non-expansive ($Df\leq1$). The above discussion provides a rigorous proof of the following sharp convergence bound.
\begin{proposition}
\label{proposition_RBM}
Let $X$ be the RBM defined by \eqref{RBM}. It has a unique stationary distribution $X_\infty$ with 
$$W(X_t,X_\infty)\leq\frac{\lambda^{t-s}}{1-\lambda^s}\cdot\frac{\E\abs{e^{cX_s}-e^{cX_0}}}{c}$$
where $t>s>0$, $c=r/\sigma^2$, and $\lambda=e^{-r^2/(2\sigma^2)}$.
\end{proposition}

It is difficult for the D\&M or D\&C methods to achieve this exact convergence rate, which equals to the drift rate ($\E_x V(Z_s)=\lambda^s V(x)$), because the downward drift is blocked by the boundary, let alone the minorization or contraction condition.

\subsection{Tandem stochastic fluid network}
To conclude this paper, we use the boundary removal technique to study a multidimensional Markov chain, which is the workload vector of a tandem stochastic fluid network (\cite{kella1992tandem}). Consider $d$ stations $s_1,...,s_d$ in series. External fluid workload only arrives at $s_1$ and is sequentially processed by $s_2,...,s_d$. Let $r_i$ be the maximal processing rate of $s_i$. The external input follows a compound renewal process where a random amount of fluid $Z$ arrives after a random length of time $T$ has passed since the last arrival. Let $\bar{X}_t$ be the remaining workload vector at time $t$, i.e., there is $\bar{X}_t^i$ amount of remaining workload in the buffer (with infinite capacity) of $s_i$. Starting from $\bar{X}_0=x\in\R_+^d$, if there is no further external input to the system, then $\bar{X}_t$ will move toward the origin along a deterministic path. We use $w(t;x)$ to denote this path, i.e., starting from $x$, without any further input, the remaining workload vector after time $t$ is $w(t;x)$. For example, let $d=2$, $r=(2,1)$, and $x=(3,0)$. Then 
\begin{equation*}
w(t;x) = \left\{
        \begin{array}{lll}
        (3,0)-(2,1)t+(0,2)t& \quad t\in[0,3/2)\\
        (0,3/2)-(0,1)(t-3/2) & \quad t\in[3/2,3)\\
        (0,0) & \quad t\in[3,\infty)
        \end{array}
    \right..
\end{equation*}
Let $T_1$ be the next arrival time and $Z_1$ be the next arrival amount. Starting from $\bar{X}_0$, we have $\bar{X}_t=w(t;\bar{X}_0)$ for $t\in[0,T_1)$ and $\bar{X}_{T_1}=w(T_1;\bar{X}_0)+(Z_1,0,...,0).$ Let $X_n$ be the remaining workload after the $n$-th arrival, i.e., $X_n=\bar{X}_{S_n}$ where $S_n=T_1+...+T_n.$ Then $X$ is a Markov chain and its random mapping representation is 
\begin{equation}
\label{def_tandem}
    f(x)=w(T;x)+\bar{Z},\;\;\bar{Z}=(Z,0,...,0).
\end{equation}
We bound the convergence under the natural stability condition
\begin{equation}
\label{stable}
r_*=\min_{i\in[d]}r_i>\E Z/\E T
\end{equation}
where $[d]=\{1,...,d\}$ and $r_*$ is the ``bottleneck'' processing rate. Before presenting the convergence bound, we find the absorbing set of $X$, i.e., $X_0\in A$ implies $X_n\in A$ for all $n\geq0$. Let $i_*=\min\{i\in[d]:r_i=r_*\}$ the (smallest) index of the bottleneck. Then the absorbing set is
$$A=\{x\in\R_+^d:x_{i_*+1}=...=x_d=0\},$$
because any station after the bottleneck that starts empty remains empty.

\begin{proposition}
\label{proposition_tandem}
Let $X$ be the Markov chain defined by \eqref{def_tandem} starting from $X_0\in A$. Under \eqref{stable}, if $\E e^{\zeta Z}<\infty$ for some $\zeta>0$, then $X$ has a unique stationary distribution $X_\infty$ with
$$W(X_n,X_\infty)\leq\frac{\lambda_*^n}{1-\lambda_*}\cdot\E\sbk{\norm{X_1-X_0}_1\cdot\frac{\exp(a_*\mathbf{1}^\top X_1)-\exp(a_*\mathbf{1}^\top X_0)}{a_*\mathbf{1}^\top X_1-a_*\mathbf{1}^\top X_0}}$$
where $\mathbf{1}$ is the all-one vector, $\norm{x}_1=\sum_{i=1}^d|x_i|$ is the $L_1$ norm, and $(a_*,\lambda_*)$ satisfies
$$\lambda_*=\E\exp(a_*(Z-r_*T))=\inf_{a\in[0,\zeta]}\E\exp(a(Z-r_*T)).$$
\end{proposition}

The exponential rate $\lambda_*$ is determined by the difference between the total input $Z$ and output $r_*T$, which means that the workload vector $X$ converges as fast as the total workload $\mathbf{1}^\top X$. Similar to Proposition \ref{proposition_RBM}, for the total workload, the optimal drift rate $\lambda_*$ ($\E_x V(X_1)\leq\lambda_* V(x)$ where $V(x)=e^{a_*x}$) is the exact rate of convergence. Since $X$ cannot converge faster than $\mathbf{1}^\top X$, for the workload vector $X$, $\lambda_*$ is also the exact rate of convergence.

The pre-multiplier is more interesting as it contains not only the total workload $\mathbf{1}^\top X$ but also $\norm{X_1-X_0}_1$, which describes how the system structure, beyond the total input and output, affects the convergence. For example, let $d=2$ and $X_0=(M,0)$ where $M$ is so large that $s_1$ cannot be depleted before the first arrival. If $r_1<r_2$, then $s_2$ is always empty, so $\norm{X_1-X_0}=\abs{\mathbf{1}^\top X_1-\mathbf{1}^\top X_0}$. If $r_1>r_2$, then the workload at $s_2$ increases while the workload at $s_1$ decreases, so $\norm{X_1-X_0}>\abs{\mathbf{1}^\top X_1-\mathbf{1}^\top X_0}$, which leads to a larger pre-multiplier. In general, the earlier the bottleneck station appears, the faster the tandem system converges. 

\subsection{Priority queues}
Interestingly, when we use the boundary removal technique to study different queueing systems, we may obtain similar convergence bounds. Consider a system with one server but $d$ queues $q_1,...,q_d$. The external input follows a $d$-dimensional compound renewal process where a random vector amount of fluid $Z$ arrives after a random length of time $T$ has passed since the last arrival ($Z_i$ arrives at $q_i$). The server operates under a priority scheme to process the workload, where queues with smaller indices have higher priorities. This means that as long as the system is not empty, the server always serves the non-empty queue with the smallest index. Let $r$ be the service rate. The stability condition is $\mathbf{1}^\top\E Z<r\E T$. Although the stability condition is independent of the priority scheme, a poor priority scheme should make the system less reliable. A reliable system can swiftly recover from unusual disturbances (quickly converge to $X_\infty$ from unusual $X_0$). The pre-multiplier in our bound can quantify how different priority schemes affect reliability.

Similar to the previous section, let $X_n$ be the remaining workload after the $n$-th arrival. By repeating the proof of Proposition \ref{proposition_tandem} verbatim, we obtain the bound in Proposition \ref{proposition_tandem} again but with
$$\lambda_*=\E\exp(a_*(\mathbf{1}^\top Z-rT))=\inf_{a\in[0,\zeta]}\E\exp(a(\mathbf{1}^\top Z-rT)).$$
The two different systems satisfy the same bound, but $\norm{X_1-X_0}_1$ in the pre-multiplier captures the structural difference between them. Now we explain why a poor priority scheme leads to a large pre-multiplier. Let $d=2$ and $X_0=(M,0)$ where $M$ is so large that the server focuses on $q_1$ before the first arrival. Then $$\norm{X_1-X_0}_1=|Z_1-rT|+Z_2=Z_1+Z_2-rT+2(Z_1-rT)_-,$$
which is larger when the busier queue is incorrectly given the lower priority ($Z_1<Z_2$).
\begin{remark}
    The above two examples can be viewed as single server queues with different inner structures, so the empty state (a single point) is the actual boundary of their state spaces. In this case, the boundary removal technique, only utilizing the all-directional contraction caused by system depletion, can lead to the optimal convergence rate. However, for general high-dimensional queueing networks (e.g., high-dimensional RBM), system depletion rarely happens, and the boundary is formed by hyperplanes. In this case, contraction happens but not simultaneously in all directions, so it cannot be captured by the local Lipschitz constant $Df(x)$. In our ongoing work, a ``directional'' CD is being developed to describe contraction in different directions.
\end{remark}

\section{From pen and paper to deep learning}
\label{section_deep}
The goal of this section is to briefly establish that our CD methodology also lends itself to the development of an automatic computational framework to bound the convergence of general state-space Markov chains. From the polynomial bound in Theorem \ref{theorem_polynomial_b} to the exponential bound in Theorem \ref{theorem_exponential}, for the first time, bounds are {\it explicitly} linked to computed functions ($U,V$ in $KV\leq V-U$). Deep learning has demonstrated superior capability in approximating functions, particularly in high-dimensional spaces. The {\it Deep Contractive Drift Calculator} (DCDC), recently introduced in \cite{qu2024deep}, is the first general-purpose, sample-based algorithm to bound the convergence of Markov chains. The DCDC approach builds upon the theoretical developments in this paper, and it highlights the appeal of having a single unified analytical condition with functions that can be parameterized ($U$ and $V$), which is the core of the CD approach that we introduce. Here, we present a summary of the DCDC algorithm.
\begin{enumerate}[label=\roman*)]
\item The contractive drift condition, an inequality by definition, is actually an equality by nature; that is, if $KV\leq V-U$ has a solution, then $KV= V-U$ also has a solution (Theorem 1 in \cite{qu2024deep}). This equality is called the contractive drift equation (CDE).
\item Inspired by the success of physics-informed neural networks (PINNs) in solving PDEs \citep{sirignano2018dgm,raissi2019physics}, DCDC solves CDEs by training neural networks and converts solutions into convergence bounds.
\item The training process is a standard application of stochastic gradient descent (SGD) where the initial location $X_0$ and the first transition $X_1=f_1(X_0)$ are repeatedly sampled. The local Lipschitz constant $Df_1(X_0)$ can be computed via automatic differentiation.
\item The effectiveness of the algorithm is illustrated by generating numerical convergence bounds for multidimensional Markov chains arising from queueing theory as well as stochastic optimization.
\end{enumerate}
The CD framework distinguishes itself from existing methods by enabling the use of deep learning for Markov chain convergence analysis. 
In the current paper, we have developed sharp convergence bounds for stylized (structured) Markov chains, increasing our confidence when applying DCDC to generate numerical convergence bounds for realistic (less structured) Markov chains. In our ongoing work, DCDC is being used to bound the convergence of the Albert and Chib’s algorithm for probit regression on real datasets (\cite{albert1993bayesian}).

\section{Proofs}
\label{section_proofs}

\subsection{Proofs for Section \ref{section_main}}
\label{proof_section_main}
\begin{lemma}
\label{lemma_change}
    Let $f:\X\rightarrow\X$ be a locally Lipschitz mapping, i.e., $Df<\infty$ everywhere. Let $\gamma:[0,1]\rightarrow\X$ be a rectifiable curve. Let $g:\X\rightarrow\R_+$ be Borel measurable function. Then
    $$L(f(\gamma);g)\leq L(\gamma;Df\cdot g\circ f).$$
\end{lemma}

\begin{proof}[Proof of Lemma \ref{lemma_change}.]
    The goal is to prove
    $$\int_0^1g(f(\gamma(t)))dL(f(\gamma|_{[0,t]}))\leq\int_0^1 g(f(\gamma(t)))Df(\gamma(t))d L(\gamma|_{[0,t]}).$$
    The two continuously increasing functions $L(f(\gamma|_{[0,t]}))$ and $L(\gamma|_{[0,t]})$ induce two Borel measures $\mu$ and $\nu$ on $[0,1]$ such that for $0\leq a<b\leq 1$,
    $$\mu((a,b])=L(f(\gamma|_{[0,b]}))-L(f(\gamma|_{[0,a]})),\;\;\nu((a,b])=L(\gamma|_{[0,b]})-L(\gamma|_{[0,a]}).$$
    Note that $\nu$ is finite because $\gamma$ is rectifiable. The first step is to show that $\mu$ is absolutely continuous with respect to $\nu$. The second step is to show that the Radon–Nikodym derivative $d\mu/d\nu$ is bounded by $Df$, the local Lipschitz constant. 
    
    To begin, we fix an $\epsilon_0>0$. For each $t\in[0,1]$, by the definition of $Df$, there exists $\eta_t>0$ such that $d(f(y),f(z))\leq(Df(\gamma(t))+\epsilon_0)d(y,z)$ for all $y,z\in B_{\eta_t}(\gamma(t))$. By the continuity of $\gamma$, there exists $\delta_t>0$ such that $\gamma|_{I(t;\delta_t)}\subset B_{\eta_t}(\gamma(t))$ where $I(t;\delta_t)=(t-\delta_t,t+\delta_t)\cap[0,1]$. These intervals form an open cover of $[0,1]$, so there exists a finite sub-cover $\{I(t_k;\delta_{t_k}),1\leq k\leq m\}$. For $0\leq a<b\leq 1$, we can insert finitely many points between $a$ and $b$ such that any pair of adjacent points belongs to one of those $m$ intervals. 
    Then we have $$\mu((a,b))=L(f(\gamma|_{(a,b)}))\leq ML(\gamma|_{(a,b)})=M\nu((a,b)),\;\;M=\max_{1\leq k\leq m}Df(\gamma(t_k))+\epsilon_0.$$ For any Borel set $B\subset[0,1]$ and $\epsilon>0$, there exists an open set $B_\epsilon$ such that $B\subset B_\epsilon$ and $\nu(B_\epsilon)\leq\nu(B)+\epsilon;$ see, for example, Theorem 1.1 of \cite{billingsley2013convergence}.
    Since every open set in $\R$ is a countable union of disjoint open intervals, we have $\mu(B_\epsilon)\leq M\nu(B_\epsilon)$.
    Then
    $$\mu(B)\leq\mu(B_\epsilon)\leq M\nu(B_\epsilon)\leq M\nu(B)+M\epsilon.$$
    By sending $\epsilon\downarrow0$, we have $\mu(B)\leq M\nu(B)$ for all Borel set $B\subset[0,1]$, which implies that $\mu$ is absolutely continuous with respect to $\nu.$
    By Theorem 5.8.8 in \cite{bogachev2007measure}, the Radon–Nikodym derivative is well-defined and satisfies
    $$\frac{d\mu}{d\nu}(t)=\lim_{\Delta t\gtz}\frac{\mu(I(t;\Delta t))}{\nu(I(t;\Delta t))}=\lim_{\Delta t\gtz}\frac{L(f(\gamma|_{I(t;\Delta t)}))}{L(\gamma|_{I(t;\Delta t)})},\;\;\nu\text{-a.e.}\;t.$$ 
    When $\Delta t\leq\delta_t$, we have $\gamma|_{I(t;\Delta t)}\subset B_{\eta_t}(\gamma(t))$ and
    \begin{align*}
        L(f(\gamma|_{I(t;\Delta t)}))=&\sup_{t_*=t_0<t_1...<t_n=t^*,\;n\geq1}\;\;\sum_{k=1}^nd(f(\gamma(t_{k-1})),f(\gamma(t_{k})))\\
        \leq&\sup_{t_*=t_0<t_1...<t_n=t^*,\;n\geq1}\;\;\sum_{k=1}^nd(\gamma(t_{k-1}),\gamma(t_{k}))(Df(\gamma(t))+\epsilon_0)\\
        =&L(\gamma|_{I(t,\Delta t)})(Df(\gamma(t))+\epsilon_0).
    \end{align*}
    where $t_*=\max(t-\Delta t,0)$ and $t^*=\min(t+\Delta t,1)$. By sending $\epsilon_0\downarrow0,$ we have 
    $$\frac{d\mu}{d\nu}(t)\leq Df(\gamma(t)),\;\;\nu\text{-a.e.}\;t.$$
\end{proof}

\begin{proof}[Proof of Theorem \ref{theorem_polynomial_b}.]
Let $m=\lceil b\rceil\geq2$.
Starting from $KV\leq V-U^{1/b}V^{1-1/b}$, we use induction to construct a set of $m$ contractive drift conditions. Suppose that we already have 
\begin{equation}
\label{kth}
    KU^{k/b}V^{1-k/b}\leq U^{k/b}V^{1-k/b}-a_kU^{(k+1)/b}V^{1-(k+1)/b}
\end{equation}
where integer $k\geq 0$ and $a_k>0$. (Clearly, $a_0=1$.) As long as $k\leq m-2$, 
by $KU\leq U$,
$$\frac{k+1}{b}=\frac{1}{b-k}+\frac{k}{b}\prs{1-\frac{1}{b-k}},\;\;1-\frac{k+1}{b}=\prs{1-\frac{k}{b}}\prs{1-\frac{1}{b-k}},$$
and H\"older's inequality,
we have
\begin{align*}
&KU^{(k+1)/b}V^{1-(k+1)/b}\\
\leq& \prs{KU}^{1/(b-k)}\prs{KU^{k/b}V^{1-k/b}}^{1-1/(b-k)}\\
\leq& U^{1/(b-k)}\prs{U^{k/b}V^{1-k/b}-a_kU^{(k+1)/b}V^{1-(k+1)/b}}^{1-1/(b-k)}\\
=&\prs{U\prs{U^{k/b}V^{1-k/b}-a_kU^{(k+1)/b}V^{1-(k+1)/b}}^{b-k-1}}^{1/(b-k)}\\
=&\prs{U^{(k+1)/b}V^{1-(k+1)/b}\prs{U^{(k+1)/b}V^{1-(k+1)/b}-a_kU^{(k+2)/b}V^{1-(k+2)/b}}^{b-k-1}}^{1/(b-k)}\\
\leq&U^{(k+1)/b}V^{1-(k+1)/b}-a_kU^{(k+2)/b}V^{1-(k+2)/b}(b-k-1)/(b-k),
\end{align*}
where we use Young's inequality
$$\frac{x^p}{p}+\frac{y^q}{q}\geq xy,\;\;x,y\geq0,\;\;p,q>1,\;\;1/p+1/q=1$$
to obtain the last line. To be specific, $$p=b-k,\;\;q=(b-k)/(b-k-1),\;\;x=\prs{U^{(k+1)/b}V^{1-(k+1)/b}}^{1/(b-k)},$$
and
$$y=\prs{\prs{U^{(k+1)/b}V^{1-(k+1)/b}-a_kU^{(k+2)/b}V^{1-(k+2)/b}}^{b-k-1}}^{1/(b-k)}$$
where the difference is non-negative because of the induction hypothesis \eqref{kth}.
Now \eqref{kth} is established for $k+1$ with 
$$a_{k+1}=a_k\cdot\frac{b-k-1}{b-k}=\dots=a_0\cdot\frac{b-1}{b}\dots\frac{b-k-1}{b-k}=\frac{b-k-1}{b}.$$ 
By induction, \eqref{kth} holds for $k=0,\dots,m-1.$
Let $F_n=f_n\circ...\circ f_1$ and
$$V_k=U^{k/b}V^{1-k/b}\cdot\prod_{l=1}^{k-1}a_l,\;\;k=0,...,m.$$
Note that the power of $V$ in $V_m$ may be negative.
Then \eqref{kth} becomes $KV_k\leq V_k-V_{k+1}$. For $n\geq1$, $x\in\X$, and $k=0,\dots,m-1$, we have
\begin{align*}
    \E DF_{n}(x)V_k(F_n(x))&\leq\E Df_n(F_{n-1}(x))V_k(f_n(F_{n-1}(x)))DF_{n-1}(x)\\
    &=\E DF_{n-1}(x)\E\sbk{Df_n(F_{n-1}(x))V_k(f_n(F_{n-1}(x)))|F_{n-1}}\\
    &\leq\E DF_{n-1}(x)(V_k(F_{n-1}(x))-V_{k+1}(F_{n-1}(x)))\\
    &=\E DF_{n-1}(x)V_k(F_{n-1}(x))-\E DF_{n-1}(x)V_{k+1}(F_{n-1}(x))\\
    &\leq...\leq V_k(x)-\sum_{l=0}^{n-1}\E DF_{l}(x)V_{k+1}(F_{l}(x)),
\end{align*}
where the first inequality is because of the submultiplicativity of the local Lipschitz constant 
$$D(g\circ h)(x)\leq Dg(h(x))Dh(x),\;\;g,h:\X\rightarrow\X,\;\;x\in\X.$$
By sending $n\gti,$ we have
\begin{equation}
\label{0step}
    \sum_{l=0}^{\infty}\E DF_{l}(x)V_{k+1}(F_{l}(x))\leq V_k(x),\;\;k=0,...,m-1.
\end{equation}
Suppose that we already have 
\begin{equation}
\label{kstep}
\sum_{n=0}^{\infty}\binom{k+n}{k}\E DF_{n}(x)V_{k+1}(F_{n}(x))\leq V(x),\;\;x\in\X,
\end{equation}
for some $k\geq 0.$ (When $k=0$, \eqref{kstep} is \eqref{0step}.) As long as $k\leq m-2$, by \eqref{0step}, we have
\begin{align*}
V(x)&\geq\sum_{n=0}^{\infty}\binom{k+n}{k}\E DF_{n}(x)V_{k+1}(F_{n}(x))\\
&\geq\sum_{n=0}^{\infty}\binom{k+n}{k}\E DF_n(x)\sum_{l=0}^{\infty}\E \sbk{D\tilde{F}_{l}(F_n(x))V_{k+2}(\tilde{F}_{l}(F_n(x)))|F_n}\\
&=\sum_{n,l=0}^{\infty}\binom{k+n}{k}\E DF_n(x) D\tilde{F}_{l}(F_n(x))V_{k+2}(\tilde{F}_{l}(F_n(x)))\\
&\geq\sum_{n,l=0}^{\infty}\binom{k+n}{k}\E DF_{n+l}(x)V_{k+2}(F_{n+l}(x))\\
&=\sum_{\bar{n}=0}^{\infty}\E DF_{\bar{n}}(x)V_{k+2}(F_{\bar{n}}(x))\sum_{l=0}^{\bar{n}}\binom{k+l}{k}\\
&=\sum_{\bar{n}=0}^{\infty}\binom{k+1+\bar{n}}{k+1}\E DF_{\bar{n}}(x)V_{k+2}(F_{\bar{n}}(x)),
\end{align*}
where $\tilde{F}_l$ is the composition of $l$ iid copies of $f$ that are independent of $F_n$. The 
combinatorial identity used in the last line can be found in \cite{gould1972combinatorial}; see (1.49) there.
Now \eqref{kstep} is established for $k+1$. By induction, \eqref{kstep} holds for $k=0,\dots,m-1$. In particular, we have
\begin{align*}
\sum_{n=0}^{\infty}\binom{m-1+n}{m-1}\E DF_{n}(x)V_{m}(F_{n}(x))\leq V(x),\;\;\sum_{n=0}^{\infty}\E DF_{n}(x)V_{1}(F_{n}(x))\leq V(x).
\end{align*}
Let $1/p=(b-1)/(m-1),\;1/q=(m-b)/(m-1).$ Again, by Young's inequality,
\begin{align*}
    &\frac{\prs{\prs{\binom{m-1+n}{m-1} DF_{n}(x)V_{m}(F_{n}(x))}^{1/p}}^p}{p}+\frac{\prs{\prs{DF_{n}(x)V_{1}(F_{n}(x))}^{1/q}}^{q}}{q}\\
    \geq&\prs{\binom{m-1+n}{m-1} DF_{n}(x)V_{m}(F_{n}(x))}^{1/p}\prs{DF_{n}(x)V_{1}(F_{n}(x))}^{1/q}\\
    =&\prs{\binom{m-1+n}{m-1} DF_{n}(x)U(F_{n}(x))^{m/b}V(F_{n}(x))^{1-m/b}\cdot\prod_{l=1}^{m-1}a_l}^{1/p}\\
    &\cdot\prs{DF_{n}(x)U(F_{n}(x))^{1/b}V(F_{n}(x))^{1-1/b}}^{1/q}\\
    =&\sbk{\binom{m-1+n}{m-1} \prod_{l=1}^{m-1}\frac{b-l}{b}}^{1/p}\\
    &\cdot DF_n(x)^{1/p+1/q}U(F_n(x))^{m/(bp)+1/(bq)}V(F_n(x))^{1/p-m/(bp)+1/q-1/(bq)}\\
    =&\sbk{\prod_{l=1}^{m-1}\frac{n+l}{m-l}\cdot\frac{b-l}{b}}^{\frac{b-1}{m-1}}\\
    &\cdot DF_n(x)U(F_n(x))^{(m(b-1)+m-b)/(b(m-1))}V(F_n(x))^{1-(m(b-1)+m-b)/(b(m-1))}\\
    =&\sbk{\prod_{l=1}^{m-1}\frac{n+l}{b}\cdot\frac{b-l}{m-l}}^{\frac{b-1}{m-1}}\cdot DF_n(x)U(F_n(x)).
\end{align*}
Let $c_n$ be the first term in the last line, which is $O(n^{b-1})$.
Then
\begin{equation}
\begin{aligned}
\label{cn}
    V(x)&=(1/p+1/q)V(x)\\
    &\geq\sum_{n=0}^{\infty}\binom{m-1+n}{m-1}\E DF_{n}(x)V_{m}(F_{n}(x))/p+\sum_{n=0}^{\infty}\E DF_{n}(x)V_{1}(F_{n}(x))/q\\
    &\geq\sum_{n=0}^{\infty}c_n\E DF_n(x)U(F_n(x)).
\end{aligned}
\end{equation}
Let $\bar{F}_n=f_1\circ...\circ f_n$ and $\bar{X}_n=\bar{F}_n(X_0)$. 
For $x,y\in\X$ and $\tilde{\gamma}\in\Gamma(x,y)$, by $\bar{F}_n(\tilde{\gamma})\in\Gamma(\bar{F}_n(x),\bar{F}_n(y))$ and Lemma \ref{lemma_change}, we have
\begin{align*}
    d_U(\bar{F}_n(x),\bar{F}_n(y))=&\inf_{\gamma\in\Gamma(\bar{F}_n(x),\bar{F}_n(y))}L(\gamma;U)\\
    \leq&L(\bar{F}_n(\tilde{\gamma});U)\\
    \leq&L(\tilde{\gamma},D\bar{F}_n\cdot U\circ\bar{F}_n)\\
    =&\int_0^1 U(\bar{F}_n(\tilde{\gamma}(t)))D\bar{F}_n(\tilde{\gamma}(t))dL(\tilde{\gamma}|_{[0,t]}),
\end{align*}
where $d_U(\bar{F}_n(x),\bar{F}_n(y))$ is measurable because it is a continuous function of two random variables. By taking expectation over $\bar{F}_n$,
\begin{equation}
\label{E_inf}
  \E d_U(\bar{F}_n(x),\bar{F}_n(y))\leq\inf_{\gamma\in\Gamma(x,y)}\int_0^1 \E U(\bar{F}_n(\gamma(t)))D\bar{F}_n(\gamma(t))dL(\gamma|_{[0,t]}).  
\end{equation}
By \eqref{cn} with $F_n$ replaced by $\bar{F}_n$ (they have the same marginal distribution),
\begin{align*}
    \sum_{n=0}^\infty c_n\E d_U(\bar{F}_n(x),\bar{F}_n(y))\leq&\inf_{\gamma\in\Gamma(x,y)}\int_0^1\sum_{n=0}^\infty c_n\E U(\bar{F}_n(\gamma(t)))D\bar{F}_n(\gamma(t))dL(\gamma|_{[0,t]})\\
    \leq&\inf_{\gamma\in\Gamma(x,y)}\int_0^1 V(\gamma(t))dL(\gamma|_{[0,t]})\\
    =&d_V(x,y).
\end{align*}
By the above inequality with $y$ replaced by $f(x)$ ($d_V(x,f(x))$ is measurable as a continuous function of a random variable), for $x\in \X,$ we have
\begin{align*}
    c_0\sum_{n=0}^\infty \E d_U(\bar{F}_n(x),\bar{F}_{n+1}(x))\leq&\sum_{n=0}^\infty c_n\E d_{U}(\bar{F}_n(x),\bar{F}_{n+1}(x))\\
    =& \sum_{n=0}^\infty c_n\E\E\sbk{d_{U}(\bar{F}_n(x),\bar{F}_{n}(f_{n+1}(x)))\Big|f_{n+1}}\\
    =& \E\sum_{n=0}^\infty c_n\E\sbk{d_{U}(\bar{F}_n(x),\bar{F}_{n}(f(x)))\Big|f}\\
    \leq& \E d_V(x,f(x)).
\end{align*}
By integrating the above inequality with respect to $X_0$, we have
$$c_0\sum_{n=0}^\infty \E d_U(\bar{X}_n,\bar{X}_{n+1})\leq\sum_{n=0}^\infty c_n\E d_U(\bar{X}_n,\bar{X}_{n+1})\leq \E d_V(X_0,X_1)<\infty.$$
This implies that $\sum_{n=0}^\infty d_U(\bar{X}_n,\bar{X}_{n+1})$ and $\sum_{n=0}^\infty d_I(\bar{X}_n,\bar{X}_{n+1})$ are finite almost surely. Since $(\X,d_I)$ is complete, wp1 there exists $\bar{X}_\infty$ such that $\lim_{n\gti}d_I(\bar{X}_n,\bar{X}_\infty)=0$. For each $n$, there exists some curve $\gamma_n$ from $\bar{X}_n$ to $\bar{X}_{n+1}$ such that $L(\gamma_n;U)<d_U(\bar{X}_n,\bar{X}_{n+1})+1/2^n$. Then $\gamma^*=\cup_{n=0}^\infty\gamma_n$ is a $d_U$-rectifiable curve from $\bar{X}_0$ to $\bar{X}_\infty$. 
As mentioned in the Preliminaries, the length function of a rectifiable curve is continuous,
so we have $\lim_{t\rightarrow1}L(\gamma^*|_{[t,1]};U)=0$,
$\lim_{n\gti}d_U(\bar{X}_n,\bar{X}_\infty)=0$,
and $$d_U(\bar{X}_n,\bar{X}_\infty)\leq\lim_{m\gti}\prs{\sum_{k=n}^{m-1} d_U(\bar{X}_k,\bar{X}_{k+1})+d_U(\bar{X}_m,\bar{X}_\infty)}=\sum_{k=n}^{\infty} d_U(\bar{X}_k,\bar{X}_{k+1}).$$
Finally,
\begin{align*}
    W_U(X_n,X_\infty)&\leq \E d_U(\bar{X}_n,\bar{X}_\infty)\\
    &\leq \sum_{k=n}^\infty\E d_U(\bar{X}_k,\bar{X}_{k+1})\\
    &\leq (1/c_n)\sum_{k=n}^\infty c_k\E d_U(\bar{X}_k,\bar{X}_{k+1})\\
    &\leq\sbk{\prod_{l=1}^{m-1}\frac{b}{n+l}\cdot\frac{m-l}{b-l}}^{\frac{b-1}{m-1}}\cdot\E d_V(X_0,X_1).
\end{align*}
Moreover, the third line implies $c_nW_U(X_n,X_\infty)\gtz$ and $W_U(X_n,X_\infty)=o(1/n^{b-1})$.
\end{proof}

\begin{proof}[Proof of Theorem \ref{theorem_semi}.]
For each integer $m\geq2,$
$$\frac{\delta V}{(\log V)^\lambda}=\frac{\delta V}{V^{1/m}}\frac{V^{1/m}}{(m\log V^{1/m})^\lambda}\geq\frac{\delta V^{1-1/m}}{m^\lambda}\cdot\inf_{x>1}\frac{x}{(\log x)^\lambda}=\frac{\delta V^{1-1/m}}{m^\lambda}\cdot\frac{e^\lambda}{\lambda^\lambda}=\frac{\kappa V^{1-1/m}}{m^\lambda}$$
where $\kappa=\delta(e/\lambda)^\lambda,$ 
so $$KV\leq V-\delta V/(\log V)^\lambda\leq V-(\kappa/m^\lambda)V^{1-1/m}.$$
By Corollary \ref{corollary_polynomial_m},
\begin{align*}
    W_I(X_n,X_\infty)\leq&\frac{m^{m\lambda}}{\kappa^m}\sbk{\prod_{k=1}^{m-1}\frac{m}{n+k}}\cdot\E d_V(X_0,X_1)\\
    =&\frac{m^{m\lambda+m}}{\kappa^mn^{m}}\frac{n}{m}\sbk{\prod_{k=1}^{m-1}\frac{1}{1+k/n}}\cdot\E d_V(X_0,X_1)\\
    \leq&\frac{m^{m\lambda+m}}{\kappa^mn^{m}}\frac{n}{m}\cdot\E d_V(X_0,X_1)\\
    =&\prs{\frac{m^{\lambda +1}}{\kappa n}}^{m}\frac{n}{m}\cdot\E d_V(X_0,X_1).
\end{align*}
To make it decay at a semi-exponential rate (e.g., $\exp(-\sqrt{n})$), we can let $m=m(n)$ increase at a certain rate that makes the expression in the parenthesis converge to some constant as $n\gti$, which suggests $m=O(n^\eta)$ where $\eta=1/(\lambda+1)$. Therefore, with $m=\lfloor cn^\eta\rfloor$ and
$$\prs{\frac{m^{\lambda +1}}{\kappa n}}^{m}\leq\prs{\frac{(cn^\eta)^{1/\eta}}{\kappa n}}^{m}=\prs{\frac{c^{1/\eta}}{\kappa }}^{m}\leq\prs{\frac{c^{1/\eta}}{\kappa }}^{cn^\eta-1}=\frac{\kappa}{c^{1/\eta}}\prs{\prs{\frac{c^{1/\eta}}{\kappa }}^c}^{n^\eta},$$
we minimizes the expression in the rightmost parenthesis
$$\log\prs{\prs{\frac{c^{1/\eta}}{\kappa }}^c}=\frac{\kappa^\eta}{\eta}\frac{c}{\kappa^\eta}\log\prs{\frac{c}{\kappa^\eta}}\geq-\frac{\kappa^\eta}{e\eta},$$
where the minimum is reached at $c=\kappa^\eta/e.$ Finally, when 
$$n\geq(2e)^{1/\eta}/\kappa\;\Rightarrow\;n^\eta\geq 2e/\kappa^\eta\;\Rightarrow\;cn^\eta\geq2\;\Rightarrow\;m\geq2,$$
we have
\begin{align*}
    W_I(X_n,X_\infty)\leq&\prs{\frac{m^{\lambda +1}}{\kappa n}}^{m}\frac{n}{m}\cdot\E d_V(X_0,X_1)\\
    \leq&\frac{\kappa}{c^{1/\eta}}\prs{e^{-\kappa^\eta/(e\eta)}}^{n^\eta}\frac{n}{2}\cdot\E d_V(X_0,X_1)\\
    =&e^{1/\eta}(n/2)e^{-n^\eta\kappa^\eta/(e\eta)}\cdot\E d_V(X_0,X_1).
\end{align*}
\end{proof}

\begin{proof}[Proof of Theorem \ref{theorem_exponential}.]
    Let $F_n=f_n\circ\dots\circ f_1$. For $n\geq1$ and $x\in\X$, we have
    \begin{align*}
        \E DF_{n}(x)V(F_n(x))&\leq\E Df_n(F_{n-1}(x))V(f_n(F_{n-1}(x)))DF_{n-1}(x)\\
        &=\E DF_{n-1}(x)\E\sbk{Df_n(F_{n-1}(x))V(f_n(F_{n-1}(x)))|F_{n-1}}\\
        &\leq r\E DF_{n-1}(x)V(F_{n-1}(x))\\
        &\leq\dots\leq r^nV(x).
    \end{align*}
    Let $\bar{F}_n=f_1\circ\dots\circ f_n$. 
    By \eqref{E_inf} with $U$ replaced by $V$, for $x,y\in\X$, we have
    \begin{align*}
        \E d_V(\bar{F}_n(x),\bar{F}_n(y))\leq&\inf_{\gamma\in\Gamma(x,y)}\int_0^1 \E V(\bar{F}_n(\gamma(t)))D\bar{F}_n(\gamma(t))dL(\gamma|_{[0,t]})\\
        \leq&\inf_{\gamma\in\Gamma(x,y)}\int_0^1 r^nV(\gamma(t))dL(\gamma|_{[0,t]})\\
        =&r^n d_V(x,y).
    \end{align*}
    By the above inequality with $y$ replaced by $f(x)$, for $x\in\X$, we have
    \begin{align*}
        \E d_V(\bar{F}_n(x),\bar{F}_{n+1}(x))=&\E\E\sbk{d_{V}(\bar{F}_n(x),\bar{F}_{n}(f_{n+1}(x)))\Big|f_{n+1}}\\
        =&\E\E\sbk{d_{V}(\bar{F}_n(x),\bar{F}_{n}(f(x)))\Big|f}\\
        \leq&r^n\E d_V(x,f(x)).
    \end{align*}
By integrating the above inequality with respect to $X_0$, we have
$$\sum_{n=0}^\infty \E d_V(\bar{X}_n,\bar{X}_{n+1})\leq \sum_{n=0}^\infty r^n\E d_V(X_0,X_1)<\infty.$$
As in the proof of Theorem \ref{theorem_polynomial_b}, wp1 there exists $\bar{X}_\infty$ such that $\lim_{n\gti}d_V(\bar{X}_n,\bar{X}_\infty)=0.$ Finally,
\begin{align*}
    W_V(X_n,X_\infty)\leq \E d_V(\bar{X}_n,\bar{X}_\infty)\leq \sum_{k=n}^\infty\E d_V(\bar{X}_k,\bar{X}_{k+1})\leq \sum_{k=n}^\infty r^k\E d_V(X_0,X_1).
\end{align*}
\end{proof}

\subsection{Proofs for Section \ref{section_expansive}}
\label{proof_section_expansive}
\begin{proof}[Proof of Proposition \ref{proposition_expansive}.]
\label{proof_proposition_expansive}
To begin, we verify $KU\leq U$ where $U(x)=x^2+1/a$. 
By symmetry, we focus on $x\geq0$.
For $x\in[0,L)$,
\begin{align*}
    KU(x)-U(x)=&\prs{1-2\gamma a}\prs{\E\sbk{x-2\gamma ax+\sqrt{2\gamma}Z}^2+1/a}-x^2-1/a\\
    =&\prs{1-2\gamma a}\prs{(1-2\gamma a)^2x^2+2\gamma+1/a}-x^2-1/a\\
    =&((1-2\gamma a)^3-1)x^2+(1-2\gamma a)2\gamma-2\gamma a(1/a)\\
    =&2\gamma\prs{1-2\gamma a-a(1/a)}\\
    \leq&-4\gamma^2 a
\end{align*}
where the last line explains why we add $1/a$ to $x^2$.
For $x\geq L$,
\begin{align*}
    KU(x)-U(x)=&\prs{1+\frac{\gamma}{4x^{3/2}}}\prs{\E\sbk{x-\frac{\gamma}{2\sqrt{x}}+\sqrt{2\gamma}Z}^2+1/a}-x^2-1/a\\
    =&\prs{1+\frac{\gamma}{4x^{3/2}}}\prs{x^2-\gamma\sqrt{x}+\frac{\gamma^2}{4x}+2\gamma+1/a}-x^2-1/a\\
    =&\frac{\gamma}{4x^{3/2}}\prs{x^2-\gamma\sqrt{x}+\frac{\gamma^2}{4x}+2\gamma+1/a}+\prs{-\gamma\sqrt{x}+\frac{\gamma^2}{4x}+2\gamma}\\
    =&\gamma\prs{\frac{\sqrt{x}}{4}-\frac{\gamma}{4x}+\frac{\gamma^2}{16x^{5/2}}+\frac{\gamma}{2x^{3/2}}+\frac{1/a}{4x^{3/2}}-\sqrt{x}+\frac{\gamma}{4x}+2}\\
    =&\gamma\prs{-\frac{3\sqrt{x}}{4}+\frac{\gamma^2}{16x^{5/2}}+\frac{2\gamma+4L^{3/2}}{4x^{3/2}}+2}\\
    \leq&\gamma\prs{-\frac{3\sqrt{L}}{4}+\frac{1}{16 L^{5/2}}+\frac{2}{4 L^{3/2}}+3}\\
    \leq&-8.99\gamma
\end{align*}
where we use $a=1/(4L^{3/2})$, $\sqrt{L}\geq16$, and $\gamma\leq1$ at the end. Now we have $KU\leq U$. Next, we compute $KV-V$ where $V(x)=x^m+M$ and $m,M\geq4$ will be determined later. For $x\geq L$,
\begin{align*}
    &KV(x)-V(x)\\
    =&\prs{1+\frac{\gamma}{4x^{3/2}}}\prs{\E\sbk{x-\frac{\gamma}{2\sqrt{x}}+\sqrt{2\gamma}Z}^{m}+M}-x^{m}-M\\
    \leq&\prs{1+\frac{\gamma}{4x^{3/2}}}\prs{
    \sbk{x-\frac{\gamma}{2\sqrt{x}}}^{m}+\binom{m}{2}x^{m-2}2\gamma+2^mx^{m-3}(2\gamma)^{\frac{3}{2}}\E Z^m+M}-x^{m}-M\\
    \leq&\prs{1+\frac{\gamma}{4x^{3/2}}}\prs{x^{m}-(\gamma/2)mx^{m-\frac{3}{2}}+\binom{m}{2}x^{m-2}2\gamma+C_mx^{m-3}\gamma^{\frac{3}{2}}+M}-x^{m}-M\\
    =&\frac{\gamma}{4x^{3/2}}\prs{x^{2}-(\gamma/2) mx^{m-3/2}+\binom{m}{2}x^{m-2}2\gamma+C_mx^{m-3}\gamma^{3/2}+M}\\
    &-(\gamma/2) mx^{m-3/2}+\binom{m}{2}x^{m-2}2\gamma+C_mx^{m-3}\gamma^{3/2}\\
    \leq&\frac{\gamma}{4}\prs{x^{m-3/2}+\binom{m}{2}x^{m-7/2}2\gamma+C_mx^{m-9/2}\gamma^{3/2}+Mx^{-3/2}}\\
    &-(\gamma/2)mx^{m-3/2}+\binom{m}{2}x^{m-2}2\gamma+C_mx^{m-3}\gamma^{3/2}\\
    \leq&\gamma \prs{x^{m-3/2}(1/4-m/2)+\binom{m}{2}x^{m-2}2+2C_mx^{m-3}\sqrt{\gamma}+\frac{M}{4x^{3/2}}}\\
    \leq&\gamma x^{m-3/2}\prs{\frac{1}{4}-\frac{m}{2}+\frac{m^2}{\sqrt{L}}+2C_m\frac{\sqrt{\gamma}}{L^{3/2}}+\frac{M}{4L^{m}}}\\
    \leq&\gamma x^{m-3/2}\prs{\frac{1}{4}-\frac{4}{2}+1+\frac{1}{16}+\frac{5}{8}}\\
    =&-(\gamma/16)x^{m-3/2}
\end{align*}
where we let $C_m=2^{m+2}\E Z^m$ to obtain the second inequality and we let $m=L^{1/4}$, $M=(5/2)L^{m}$, $\sqrt{L}\geq16$, $2C_m\sqrt{\gamma}\leq256$ to obtain the last inequality. Here $C_m$ and $\tilde{C}_m$ are constants that only depend on $m$. With $b=(2/3)(m-2)$, we have
\begin{align*}
    U(x)^{1/b}V(x)^{1-1/b}=&\prs{x^2+1/a}^{1/b}\prs{x^{m}+M}^{1-1/b}\\
    =&x^{m-3/2}\prs{1+\frac{4L^{3/2}}{x^2}}^{1/b}\prs{1+\frac{(5/2)L^{m}}{x^{m}}}^{1-1/b}\\
    \leq&x^{m-3/2}\prs{1+\frac{4}{16}}^{1/b}\prs{1+\frac{5}{2}}^{1-1/b}\\
    \leq&(7/2)x^{m-3/2}
\end{align*}
Now we have $KV\leq V-(\gamma/56) U^{1/b}V^{1-1/b}$ for $x\geq L$. For $x\in[0,L),$
\begin{align*}
    &KV(x)-V(x)\\
    =&\prs{1-2\gamma a}\prs{\E\sbk{x-2\gamma ax+\sqrt{2\gamma}Z}^{m}+M}-x^{m}-M\\
    =&\prs{1-2\gamma a}\E\sbk{x-2\gamma ax+\sqrt{2\gamma}Z}^{m}-x^{m}-2\gamma aM\\
    =&\sum_{k=0}^{m}\binom{m}{k}(1-2\gamma a)^{m-k+1}x^{m-k}(2\gamma)^{k/2}\E Z^k-x^{m}-2\gamma aM\\
    =&((1-2\gamma a)^{m+1}-1)x^{m}+\sum_{k=2}^{m}\binom{m}{k}(1-2\gamma a)^{m-k+1}x^{m-k}(2\gamma)^{k/2}\E Z^k-2\gamma aM\\
    \leq&\binom{m}{2}x^{m-2}2\gamma+2^m(1+x^{m-3})(2\gamma)^{3/2}-2\gamma aM\\
    \leq&-2\gamma\prs{aM-(m^2/2)x^{m-2}-\bar{C}_m\sqrt{\gamma}L^{m-3}}\\
    \leq&-2\gamma\prs{\frac{(5/2)L^{m}}{4L^{3/2}}-\frac{\sqrt{L}L^{m-2}}{2}-\frac{L^{m-3/2}}{16}}\\
    =&-(\gamma/8) L^{m-3/2}
\end{align*}
where we let $\bar{C}_m=2^{m+3/2}\E Z^m$ to obtain the second inequality and we let 
$\bar{C}_m\sqrt{\gamma}L^{m-3}\leq L^{m-3/2}/16$ to obtain the last inequality. In addition, we have
\begin{align*}
    U(x)^{1/b}V(x)^{1-1/b}=&\prs{x^2+1/a}^{1/b}\prs{x^{m}+M}^{1-1/b}\\
    \leq&\prs{L^2+1/a}^{1/b}\prs{L^{m}+M}^{1-1/b}\\
    =&L^{m-3/2}\prs{1+\frac{4L^{3/2}}{L^2}}^{1/b}\prs{1+\frac{(5/2)L^{m}}{L^{m}}}^{1-1/b}\\
    \leq&(7/2)L^{m-3/2}.
\end{align*}
Now we have $KV\leq V-(\gamma/28) U^{1/b}V^{1-1/b}$ for $x\in[0,L).$ Finally, we have $KU\leq U$ and $KV\leq V-(\gamma/56) U^{1/b}V^{1-1/b}$ hold everywhere. By Theorem \ref{theorem_polynomial_b},
\begin{align*}(\gamma/56)^b4L^{3/2}W(X_n,X_\infty)\leq&(\gamma/56)^bW_U(X_n,X_\infty)\\
\leq&\sbk{\prod_{k=1}^{\lceil b\rceil-1}\frac{b}{n+k}\cdot\frac{\lceil b\rceil-k}{b-k}}^{\frac{b-1}{\lceil b\rceil-1}}\cdot\E d_V(X_0,X_1)\\
=&\sbk{\prod_{k=1}^{\lceil b\rceil-1}\frac{b}{n+k}\cdot\frac{\lceil b\rceil-k}{b-k}}^{\frac{b-1}{\lceil b\rceil-1}}\cdot\E\sbk{\int_{X_0\wedge X_1}^{X_0\vee X_1}V(x)dx}
\end{align*}
where $b=(2/3)(L^{1/4}-2)$, $V(x)=x^{L^{1/4}}+(5/2)L^{L^{1/4}}$, and $U(x)=x^2+4L^{3/2}.$
We can let $\sqrt{\gamma}\leq 2^{13/2-m}/\E Z^m$ to make sure $2C_m\sqrt{\gamma}\leq256$ and $\bar{C}_m\sqrt{\gamma}L^{m-3}\leq L^{m-3/2}/16$.
\end{proof}

\subsection{Proofs for Section \ref{section_sgd}}
\label{proof_section_sgd}
\begin{proof}[Proof of Proposition \ref{proposition_sgd_1}.]
Let $V(x)=\delta(1-|x|)_++1$ where $\delta>0$ will be determined later.
By symmetry, we focus on $x\geq0$.
For $x\geq1$,
\begin{align*}
    KV(x)-V(x)&=\E\prs{1-\alpha}\prs{\delta(1-\abs{x-\alpha(x+Z)})_++1}-1\\
    &=-\alpha+\E\prs{1-\alpha}\delta\prs{1-\abs{(1-\alpha)x-\alpha Z}}_+\\
    &\leq-\alpha+\E\prs{1-\alpha}\delta\prs{1-(1-\alpha)x+\alpha Z}_+\\
    &\leq-\alpha+\E\delta\prs{1-(1-\alpha)+\alpha Z}_+\\
    &\leq-\alpha\prs{1-\delta\E\prs{1+Z}_+}\\
    &\leq-\alpha\prs{1-\delta\E\prs{1+|Z|}}\\
    &\leq-(3/4)\alpha
\end{align*}
where we let $\delta\leq(1/4)/\E(1+|Z|)$ to obtain the last inequality.
For $x\in[0,1)$,
\begin{align*}
    &KV(x)-V(x)\\
    =&\E\prs{1-\alpha (m-1)x^{m-2}}\prs{\delta(1-\abs{x-\alpha(x^{m-1}+Z)})_++1}-(\delta(1-x)+1)\\
    \leq&\delta\E\prs{(1-\abs{x-\alpha(x^{m-1}+Z)})_+-(1-x)}-\alpha (m-1)x^{m-2}.
\end{align*}
For the first term, 
\begin{align*}
    &\E(1-\abs{x-\alpha(x^{m-1}+Z)})_+-(1-x)\\
    =&\E(1-\abs{\alpha Z})_+-1+\E(1-\abs{x-\alpha(x^{m-1}+Z)})_+-\E(1-\abs{\alpha Z})_++x\\
    \leq&\E(1-\abs{\alpha Z})_+-1+\E\abs{(1-\abs{x-\alpha(x^{m-1}+Z)})-(1-\abs{-\alpha Z})}+x\\
    \leq&\E(1-\abs{\alpha Z})_+-1+\E\abs{-\prs{x-\alpha(x^{m-1}+Z)}+\prs{-\alpha Z}}+x\\
    \leq&-\bar{\alpha}+2x
\end{align*}
where we let $\alpha\in(0,1)$ and $\bar{\alpha}=1-\E(1-\abs{\alpha Z})_+=\Theta(\alpha)$ to obtain the last inequality. We need to choose $\delta$ to make $$KV(x)-V(x)\leq\delta(-\bar{\alpha}+2x)-\alpha(m-1)x^{m-2},\;\;x\in[0,1)$$ uniformly negative. When $m=3$ and $\delta=\alpha$, the above expression is $-\bar{\alpha}\delta.$ When $m>3$, the above expression reaches its maximum at $x_*=(2\delta/(\alpha(m-1)(m-2)))^{1/(m-3)}$ and the maximum is bounded by $\delta\prs{-\bar{\alpha}+2x_*}$, so we let $x_*=\bar{\alpha}/4$ to get $-\delta\bar{\alpha}/2$ where
$$\delta=(\bar{\alpha}/4)^{m-3}\alpha(m-1)(m-2)/2\leq 2\alpha.$$
Now we have $KV\leq V-(3/4)\alpha$ in $[1,\infty)$ and $KV\leq V-\delta\bar{\alpha}/2$ in $[0,1).$ Since $\bar{\alpha}<1$,
$$\frac{\delta\bar{\alpha}/2}{(3/4)\alpha}=\frac{(\bar{\alpha}/4)^{m-2}\alpha(m-1)(m-2)}{(3/4)\alpha}<(4/3)(1/4)^{m-2}(m-1)(m-2)<1,\;\;m\geq3.$$
Now we have $$KV\leq V-\delta\bar{\alpha}/2\leq V-(\delta\bar{\alpha}/2)(V/(1+\delta))=rV,\;\;r=1-(\delta\bar{\alpha})/(2(1+\delta))$$
everywhere. By Theorem \ref{theorem_exponential},
\begin{align*}
    W(X_n,X_\infty)\leq W_V(X_n,X_\infty)\leq\sbk{r^n/(1-r)}\cdot\E d_V(X_0,X_1).
\end{align*}
Since $V\leq1+\delta<2$, the second term $\E d_V(X_0,X_1)$ is bounded by $2\alpha\E\abs{h'(X_0)+Z_1}$. For the first term,
\begin{align*}
    \frac{r^n}{1-r}=&\prs{1-\frac{\delta\bar{\alpha}}{2(1+\delta)}}^n\frac{2(1+\delta)}{\delta\bar{\alpha}}\\
    \leq&\prs{1-\frac{\delta\bar{\alpha}}{4}}^n\frac{4}{\delta\bar{\alpha}}\\
    =&\prs{1-(\bar{\alpha}/4)^{m-2}\alpha(m-1)(m-2)/2}^n(4/(\delta\bar{\alpha}))\\
    \leq&\prs{1-(\bar{\alpha}/4)^{m-2}\alpha}^n(1/((\bar{\alpha}/4)^{m-2}\alpha)).
\end{align*}
With $\tilde{\alpha}=\bar{\alpha}/4$, we have
$$W(X_n,X_\infty)\leq(2/\tilde{\alpha}^{m-2})\cdot\prs{1-\tilde{\alpha}^{m-2}\alpha}^n\cdot \E\abs{h'(X_0)+Z_1}.$$
We can let $\alpha\leq(1/8)/\E(1+|Z|)$ to make sure $\delta\leq(1/4)/\E(1+|Z|)$ because $\delta\leq2\alpha$.
\end{proof}

\begin{proof}[Proof of Proposition \ref{proposition_sgd_2}.]
The integrability condition $\E d_V(X_0,X_1)<\infty$ in Theorem \ref{theorem_polynomial_b} suggests us to consider $V(x)=|x|^{\gamma-1}+M$ with $M>1$ because $X_1-X_0=\alpha h'(X_0)+Z_1$, $\E |Z_1|^\gamma<\infty$, and
$$d_V(X_0,X_1)=\int_{X_0\wedge X_1}^{X_0\vee X_1}(|x|^{\gamma-1}+M)dx\leq c(|Z_1|^{\gamma}+|X_0|^{\gamma}+1),\;\;c>0.$$
By symmetry, we focus on $x\geq0$.
For $x\geq1$,
\begin{align*}
    &KV(x)-V(x)\\
    =&\E(1-\alpha(\beta-1)x^{\beta-2})(|(x-\alpha(x^{\beta-1}+Z)|^{\gamma-1}+M)-(x^{\gamma-1}+M)\\
    =&\E(|x-\alpha(x^{\beta-1}+Z)|^{\gamma-1}-x^{\gamma-1})-\E\alpha(\beta-1)x^{\beta-2}(|x-\alpha(x^{\beta-1}+Z)|^{\gamma-1}+M).
\end{align*}
Since both of the two terms above are $O(x^{\gamma+\beta-3})$ as $x\gti$, we should choose $b=(\gamma-1)/(2-\beta)$ in $KV\leq V-\delta V^{1-1/b}$ ($\delta$ will be determined later) as $(x^{\gamma-1})^{1-1/b}=x^{\gamma-1-(2-\beta)}=x^{\gamma+\beta-3}$. In fact, the first term is enough to establish a CD, and it satisfies
\begin{align*}
    \frac{x^{\gamma-1}-\E|x-\alpha(x^{\beta-1}+Z)|^{\gamma-1}}{(x^{\gamma-1}+M)^{1-1/b}}
    \geq\frac{x^{\gamma-1}-\E|x-\alpha(x^{\beta-1}+Z)|^{\gamma-1}}{x^{\gamma+\beta-3}}\frac{1}{(1+M)^{1-1/b}}.
\end{align*}
When $x-\alpha(x^{\beta-1}+Z)<0,$
\begin{align*}
    &\frac{\E\prs{x^{\gamma-1}-|x-\alpha(x^{\beta-1}+Z)|^{\gamma-1}}I(x-\alpha(x^{\beta-1}+Z)<0)}{x^{\gamma+\beta-3}}\\
    \geq&-\frac{\E\prs{|x-\alpha(x^{\beta-1}+Z)|^{\gamma-1}-x^{\gamma-1}}I(x-\alpha(x^{\beta-1}+Z)<-x)}{x^{\gamma+\beta-3}}\\
    \geq&-\E\prs{|x-\alpha(x^{\beta-1}+Z)|^{\gamma-1}-x^{\gamma-1}}I(x+(x-\alpha x^{\beta-1})<\alpha Z)\\
    \geq&-\E(\alpha Z-(x-\alpha x^{\beta-1}))^{\gamma-1}I(x+(x-\alpha x^{\beta-1})<\alpha Z)\\
    \geq&-\E(\alpha Z)^{\gamma-1}I(x+(x-\alpha x^{\beta-1})<\alpha Z)\\
    \geq&-\E(\alpha Z)^{\gamma-1}I(1<\alpha Z)\\
    \geq&-\alpha^{\gamma}\E|Z|^\gamma.
\end{align*}
When $x-\alpha(x^{\beta-1}+Z)\geq0,$
\begin{align*}
    &\frac{\E\prs{x^{\gamma-1}-(x-\alpha(x^{\beta-1}+Z))^{\gamma-1}}I(x-\alpha(x^{\beta-1}+Z)\geq0)}{x^{\gamma+\beta-3}}\\
    =&\frac{\E\prs{1-(1-\alpha(x^{\beta-2}+Z/x))^{\gamma-1}}I(1-\alpha(x^{\beta-2}+Z/x)\geq0)}{x^{\beta-2}}\\
    \geq&\frac{\E\prs{\gamma-1)\alpha(x^{\beta-2}+Z/x}I(\alpha Z\leq x-\alpha x^{\beta-1})}{x^{\beta-2}}\\
    =&(\gamma-1)\alpha\prs{P(\alpha Z\leq x-\alpha x^{\beta-1})-\frac{\E (-Z)I(\alpha Z\leq x-\alpha x^{\beta-1} )}{x^{\beta-1}}}\\
    \geq&(\gamma-1)\alpha\prs{P(\alpha Z\leq 1-\alpha)-\sup_{\bar{x}\geq1}\sbk{\E (-Z)I(\alpha Z\leq \bar{x}-\alpha \bar{x}^{\beta-1})}}
\end{align*}
where the first inequality is because $(1-y)^a\leq 1-ay$ where $y<1,a\in(0,1)$ and the second inequality is because $x-\alpha x^{\beta-1}$ is increasing in $[1,\infty)$ and $\E Z=0$ implies $\E ZI(Z\leq\cdot)\leq0.$ Now for $x\geq 1$ we have
\begin{align*}
    &(1+M)^{1-1/b}\cdot\frac{V(x)-KV(x)}{V(x)^{1-1/b}}\\
    \geq&(\gamma-1)\alpha P(\alpha Z\leq 1-\alpha)-(\gamma-1)\alpha\sup_{z\geq(1-\alpha)/\alpha}\sbk{\E (-Z)I(Z\leq z)}-\alpha^\gamma\E|Z|^\gamma.
\end{align*}
Note that the positive term is $\Theta(\alpha)$ while the two negative terms are $o(\alpha)$, so when $\alpha$ is small enough the above expression is larger than $(\gamma-1)\alpha/2.$ For $x\in[0,1),$
\begin{align*}
    &(1+M)^{1-1/b}\cdot\frac{V(x)-KV(x)}{V(x)^{1-1/b}}\\
    =&(1+M)^{1-1/b}\cdot\frac{(x^{\gamma-1}+M)-\E(1-\alpha)(|x-\alpha(x+Z)|^{\gamma-1}+M)}{(x^{\gamma-1}+M)^{1-1/b}}\\
    \geq&\alpha M+x^{\gamma-1}-\E(1-\alpha)|x-\alpha(x+Z)|^{\gamma-1}\\
    \geq&\alpha M-\E(1+|Z|)^{\gamma-1}\\
    =&(\gamma-1)\alpha/2
\end{align*}
where we let $M=(\E(1+|Z|)^{\gamma-1}+(\gamma-1)\alpha/2)/\alpha$ to obtain the last equality. Now we have 
$$KV\leq V-(1+M)^{1-1/b}(\gamma-1)(\alpha/2)V^{1-1/b},\;\;b=(\gamma-1)/(2-\beta)$$
everywhere. By Theorem \ref{theorem_polynomial_b},
\begin{align*}
&(1+M)^{b-1}(\gamma-1)^b(\alpha/2)^bW(X_n,X_\infty)\\
\leq&\sbk{\prod_{k=1}^{\lceil b\rceil-1}\frac{b}{n+k}\cdot\frac{\lceil b\rceil-k}{b-k}}^{\frac{b-1}{\lceil b\rceil-1}}\cdot\E\sbk{\int_{X_0\wedge X_1}^{X_0\vee X_1}V(x)dx}.
\end{align*}
\end{proof}

\subsection{Proofs for Section \ref{section_lmt}}
\label{proof_section_lmt}
\begin{proof}[Proof of Proposition \ref{proposition_exact}.]
Let $V_M(x)=(x+M)^m$ where $M>1$ will be determined later. 
When $m=1$, it corresponds to the standard large M technique, which has been discussed at the beginning of Section \ref{section_lmt}. Now we focus on $m\geq2.$
An obvious but useful fact is that $f(x)=(x+Z)_+=0$ when $1-Df(x)=I(x+Z<0)=1$.
For $x\geq0$,
\begin{align*}
    &\frac{V_M(x)-\E Df(x)V_M(f(x))}{V_M(x)^{1-1/m}}\\
    =&\frac{\E [1-Df(x)]V_M(f(x))+\E[V_M(x)-V_M(f(x))]}{V_M(x)^{1-1/m}}\\
    =&\frac{\E I(x+Z<0)(f(x)+M)^m+\E[(x+M)^m-(f(x)+M)^m]}{(x+M)^{m-1}}\\
    =&\frac{P(x+Z<0)M^m+\E[(x+M)^m-(f(x)+M)^m]}{(x+M)^{m-1}}\\
    =&\frac{1}{(x+M)^{m-1}}\prs{P(x+Z<0)M^m-\E\sum_{k=1}^m\binom{m}{k}(f(x)-x)^k(x+M)^{m-k}}\\
    \geq&\frac{P(x+Z<0)M^m}{(x+M)^{m-1}}+m\E(x-f(x))-\frac{1}{x+M}\E\sum_{k=2}^m\binom{m}{k}|Z|^k(x+M)^{2-k}\\
    \geq&\frac{P(x+Z<0)M^m}{(x+M)^{m-1}}+m\E(x-(x+Z)_+))-\frac{\E(1+|Z|)^m}{x+M}.\\
\end{align*}
Note that the second term above is continuous and converges $m\delta=m(-\E Z)>0$ as $x\gti$. Moreover, the limit cannot be reached until $P(x+Z<0)=0$. Therefore, there exists $\bar{x}$ with $P(\bar{x}+Z<0)>0$ such that the second term above is larger than $m(-\E Z)/2$ for all $x\geq\bar{x}.$ At $\bar{x}$, if we choose $M$ such that the third term above is larger than $m(-\E Z)/4$, then the above expression (the sum of three terms) is larger than $m(-\E Z)/4$ for all $x\geq\bar{x}.$ For $x\in[0,\bar{x})$, since $P(\bar{x}+Z<0)>0$, we can increase $M$ until the sum of the first term and the third term above is larger than $m(-\E Z)/4$ for all $x\in[0,\bar{x}).$ Now we have $KV_M\leq V_M-(m(-\E Z)/4)V_M^{1-1/m}$. By Theorem \ref{theorem_polynomial_b}, $W(X_n,X_\infty)=o(1/n^{m-1})$.

Next, we show that this polynomial rate is exact. By stochastic monotonicity and Spitzer's identity (\cite{spitzer1956combinatorial}), 
$$W(X_n,X_\infty)=\E X_\infty-\E X_n=\sum_{k=n+1}^\infty\E (S_k)_+/k$$
where $S_k=\sum_{l=1}^k Z_l.$ Suppose that there exists $a,b>0$ such that for all $n\geq1$
\begin{equation}
\label{sn lowerbound}
    \E (S_n)_+\geq an^2P(Z>b(n-1)),
\end{equation}
which will be proved later.
If there exists $\epsilon>0$ such that
$$O(n^{-(m-1+2\epsilon)})=W(X_n,X_\infty)\geq\sum_{k=n+1}^\infty akP(Z>b(k-1))\geq a\int_n^\infty xP(Z>bx)dx,$$
then 
$$\int_0^\infty y^{m-2+\epsilon}\int_y^\infty xP(Z>x)dxdy<\infty$$
as $\int_0^1 y^{m-2+\epsilon}dy<\infty$ ($m\geq1$) and $\int_1^\infty y^{m-2+\epsilon-(m-1+2\epsilon)}dy=\int_1^\infty y^{-1-\epsilon}dy<\infty$. However,
\begin{align*}
     \int_0^\infty xP(Z>x)\int_0^x y^{m-2+\epsilon}dydx=&\int_0^\infty xP(Z>x)\frac{x^{m-1+\epsilon}}{m-1+\epsilon}dx\\
    =&\int_0^\infty \frac{P(Z_+^{m+1+\epsilon}>x^{m+1+\epsilon})}{(m-1+\epsilon)(m+1+\epsilon)}dx^{m+1+\epsilon}\\
    =&\frac{\E Z_+^{m+1+\epsilon}}{(m-1+\epsilon)(m+1+\epsilon)}\\
    =&\infty
\end{align*}
leads to a contradiction, so for any $\epsilon>0$, $n^{m-1+2\epsilon}W(X_n,X_\infty)$ must be unbounded. 

Now we prove \eqref{sn lowerbound}. For set $A\subset\{1,...,n\}$, let $S_n^{-A}=\sum_{k\in A^c}Z_k$. Recall that $\delta=-\E Z>0.$
Note that $S_n$ is larger than $x$ when one $Z_l$ is larger than $2(n-1)\delta$ and the sum of the rest is larger than $x-2(n-1)\delta$.
Let $Z_{\{i,j\}}=Z_i\wedge Z_j=\min(Z_i,Z_j)$. By Bonferroni's inequality (\cite{bonferroni1936teoria}),
\begin{align*}
    &\E (S_n)_+\\
    =&\int_0^\infty P(S_n>x)dx\\
    \geq&\int_0^\infty \binom{n}{1}P(Z_1>2(n-1)\delta,\;S_n^{-\{1\}}>x-2(n-1)\delta)dx\\
    &-\int_0^\infty \binom{n}{2}P(Z_i>2(n-1)\delta,\;S_n^{-\{i\}}>x-2(n-1)\delta,\;i=1,2)dx\\
    =&\int_0^\infty \binom{n}{1}P(Z>2(n-1)\delta)P(S_n^{-\{1\}}>x-2(n-1)\delta)dx\\
    &-\int_0^\infty \binom{n}{2}P(S_n^{-\{1,2\}}+Z_{\{1,2\}}>x-2(n-1)\delta,\;Z_{\{1,2\}}>2(n-1)\delta)dx\\
    =&\binom{n}{1}P(Z>2(n-1)\delta)\E(S_{n-1}+2(n-1)\delta)_+\\
    &-\binom{n}{2}P(Z>2(n-1)\delta)^2\E\sbk{(S_{n-2}+Z_{\{n-1,n\}}+2(n-1)\delta)_+|Z_{\{n-1,n\}}>2(n-1)\delta}.
\end{align*}
Since $S_n/n\stackrel{L^1}{\rightarrow}-\delta$ and $x_+$ is Lipschitz, for the first term above, we have
$$\binom{n}{1}P(Z>2(n-1)\delta)\E(S_{n-1}+2(n-1)\delta)_+\sim \delta n^2P(Z>2(n-1)\delta)$$
where $a_n\sim b_n$ means that $a_n/b_n\rightarrow1$ as $n\gti$, so the first term satisfies \eqref{sn lowerbound}. For the second term,
\begin{align*}
    &\binom{n}{2}P(Z>2(n-1)\delta)^2\E\sbk{(S_{n-2}+Z_{\{n-1,n\}}+2(n-1)\delta)_+|Z_{\{n-1,n\}}>2(n-1)\delta}\\
    \leq&\binom{n}{2}P(Z>2(n-1)\delta)^2\E\sbk{(S_{n-2}+2(n-1)\delta)_++\frac{Z_{n-1}+Z_n}{2}\Bigg|Z_{\{n-1,n\}}>2(n-1)\delta}\\
    =&\binom{n}{2}P(Z>2(n-1)\delta)^2\prs{\E(S_{n-2}+2(n-1)\delta)_++\E\sbk{Z| Z>2(n-1)\delta}}\\
    \sim&(n^2/2)P(Z>2(n-1)\delta)\prs{P(Z>2(n-1)\delta)n\delta+\E ZI(Z>2(n-1)\delta)}.
\end{align*}
Since $Z$ is integrable, both terms in the parenthesis vanish as $n\gti$. Finally,
$$\E(S_n)_+\geq\delta n^2P(Z>2(n-1)\delta)(1-o(1)),$$
so it satisfies \eqref{sn lowerbound}.
\end{proof}

\begin{proof}[Proof of Proposition \ref{proposition_uniform}.]
Let $Y^\delta=Y-\delta.$
Let $V_M(x)=|x+M|^m-M^m+c$ where $M\geq b$ and $c\in(0,M^m)$ will be determined later. 
For $x\geq0,$
\begin{align*}
    &\frac{\E Df^\delta(x)V_M(f^\delta(x))-V_M(x)}{V_M(x)^{1-1/m}}\\
    =&\frac{\E I(x+Y^\delta\geq 0)(\abs{(x+Y^\delta)_++M}^m-M^m+c)-(\abs{x+M}^m-M^m+c)}{(\abs{x+M}^m-M^m+c)^{1-1/m}}\\
    =&\frac{\E (1-I(x+Y^\delta<0))(\abs{x+Y^\delta+M}^m-M^m+c)-(\abs{x+M}^m-M^m+c)}{(\abs{x+M}^m-M^m+c)^{1-1/m}}\\
    =&\frac{\E \abs{x+Y^\delta+M}^m-\abs{x+M}^m-\E I(x+Y^\delta<0)(\abs{x+Y^\delta+M}^m-M^m+c)}{(\abs{x+M}^m-M^m+c)^{1-1/m}}\\
    \leq&\frac{\E \abs{x+Y^\delta+M}^m-\abs{x+M}^m-\E I(x+Y^\delta<0)(\abs{x+Y^\delta+M}^m-M^m+c)}{(x+M)^{m-1}}\\
    =&\E\abs{\sum_{k=0}^m\binom{m}{k}(Y^{\delta})^k(x+M)^{1-k}}-(x+M)\\
    &-\frac{\E I(x+Y^\delta<0)(\abs{x+Y^\delta+M}^m-M^m+c)}{(x+M)^{m-1}}\\
    \leq&\E\abs{mY^\delta+(x+M)}-(x+M)+\sum_{k=2}^m\binom{m}{k}\E|Y^{\delta}|^k(x+M)^{1-k}\\
    &-\frac{P(x+Y^\delta<0)}{(x+M)^{m-1}}\cdot\E\sbk{\abs{x+Y^\delta+M}^m-M^m+c\Big|x+Y^\delta<0}\\
    \leq&\E(mY^\delta+(x+M))+2\E(mY^\delta+(x+M))^--(x+M)+\frac{\E\prs{1+|Y^\delta|}^m}{x+M}\\
    &-\frac{P(x+Y^\delta<0)}{(x+M)^{m-1}}\cdot\prs{\abs{\E\sbk{x+Y^\delta\Big|x+Y^\delta\leq0}+M}^m-M^m+c}\\
    \leq&-m\delta+2P(mY^\delta+(x+M)\leq 0)\E\sbk{-mY^\delta-(x+M)\Big|mY^\delta+(x+M)\leq 0}\\
    &+\E\prs{2+|Y|}^m/M-\frac{P(x+Y^\delta<0)}{(x+M)^{m-1}}\cdot\prs{\abs{M-b}^m-M^m+c}\\
    \leq&-m\delta+2 P\prs{Y_-+1\geq(x+M)/m-\delta+1} mb\\
    &+\E\prs{2+|Y|}^m/M-\frac{P(x+Y^\delta<0)}{(x+M)^{m-1}}\cdot\prs{c-\sum_{k=1}^m\binom{m}{k}b^kM^{m-k}}\\
    \leq&-m\delta+2\frac{\E (1+Y^-)}{M/m} mb+\E\prs{2+|Y|}^m/M-\frac{P(x+Y^\delta<0)}{(x+M)^{m-1}}\cdot\prs{c-M^{m-1}(1+b)^m},
\end{align*}
where the first inequality is because of $c<M^m$, the third inequality is because of Jensen's inequality, the fourth and fifth inequalities are because of \eqref{minus_b}, and the last inequality is because of Markov's inequality.
We choose $c=M^{m-1}(1+b)^m$ to eliminate the last term above.
Then we choose $M=4\E(2+|Y|)^m(1+b)^m/\delta$ to make sure that the second term above is less than $m\delta/2$, the third term above is less than $m\delta/4$, and $c<M^m.$
Now we have $KV_M\leq V_M-(m\delta/4)V_M^{1-1/m}$. By Corollary \ref{corollary_polynomial_m},
\begin{align*}
 W(X^\delta_n,X^\delta_\infty)&\leq\frac{1}{(m\delta/4)^m}\cdot\sbk{\prod_{k=1}^{m-1}\frac{m}{n+k}}\cdot\E\sbk{\int_0^{Y^\delta_+}\sbk{(x+M)^m-M^m+c}dx}
\end{align*}
where
\begin{align*}
    \E\sbk{\int_0^{Y^\delta_+}\sbk{(x+M)^m-M^m+c}dx}\leq&\E\sbk{\frac{(Y_++M)^{m+1}-M^{m+1}}{m+1}-M^mY_++cY_+}\\
    =&\E\sbk{\frac{1}{m+1}\sum_{k=2}^{m+1}\binom{m+1}{k}Y_+^kM^{m+1-k}+cY_+}\\
    \leq&\E\sbk{\frac{M^{m-1}}{m+1}\prs{1+Y_+}^{m+1}+cY_+}.
\end{align*}
For the scaled process,
\begin{align*}
    W(\delta X^\delta_{n/\delta^2},\delta X^\delta_\infty)&\leq\frac{4/m}{(m\delta/4)^{m-1}}\cdot\sbk{\prod_{k=1}^{m-1}\frac{m}{n/\delta^2+k}}\cdot\E\sbk{\frac{M^{m-1}}{m+1}\prs{1+Y_+}^{m+1}+cY_+}\\
    &=\frac{4}{m}\sbk{\prod_{k=1}^{m-1}\frac{M/(\delta/4)}{n/\delta^2+k}}\cdot\E\sbk{\frac{\prs{1+Y_+}^{m+1}}{m+1}+(1+b)^mY_+}\\
    &\leq\frac{4}{m}\sbk{\frac{16\E(2+|Y|)^m(1+b)^m}{n}}^{m-1}\E\sbk{\frac{\prs{1+Y_+}^{m+1}}{m+1}+(1+b)^mY_+}.
\end{align*}
\end{proof}

\subsection{Proofs for Section \ref{section_brt}}
\label{proof_section_brt}
\begin{proof}[Proof of Proposition \ref{proposition_tandem}.]
Recall that the random mapping representation is $f(x)=w(x;T)+\bar{Z}$. To begin, we argue that it is non-expansive ($Df\leq1$) with respect to the $L^1$ distance $\norm{x-y}_1=\sum_{i=1}^d\abs{x_i-y_i}$.
Starting from $x,y\in\R_+^d$ that are close to each other, we have $w_i(t;x)-w_i(t;y)=x_i-y_i$ until $s_i$ is empty. After $s_i$ is empty, $w_i(t;x)-w_i(t;y)=0$ but $x_i-y_i$ is added to $w_j(t;x)-w_j(t;y)$ where $j>i$ is the index of the next non-empty station. If no such $j$ exists, then $x_i-y_i$ simply disappears when $s_i$ becomes empty. Essentially, differences at different stations merge and eventually vanish, so
$$\norm{f(x)-f(y)}_1=\norm{w(t;x)-w(t;y)}_1=\sum_{i=1}^d\abs{w_i(t;x)-w_i(t;y)}$$
never increases, and hence $Df\leq1$. Let $w_*(t;x)$ be the extension of $w(t;x)$ beyond the origin, i.e., when $w(t;x)$ stops at the origin, $w_*(t;x)$ keeps moving without changing direction. For example, if $w(\tau;x)=0$ and $w(\tau-t;0)=(\tau-t)v$ as $t\uparrow\tau$ where $v\in\R_+^d$, then $w_*(t;x)=(\tau-t)v$ for all $t\geq\tau.$ Next, we argue that the Lipschitz constant is
$$Df(x)=I(w_*(T;x)\geq0).$$
When $w_*(T;x)<0$, $w(T;\cdot)$ maps a small neighborhood of $x$ to the origin, so $Df(x)=0.$
Recall that $A$ is the absorbing set of $X$ where all stations after the bottleneck remain empty. Starting from $x\in A$, the total workload $\mathbf{1}^\top w(t;x)$ decreases at rate $r_*$ until it hits the origin.
Moreover, $\mathbf{1}^\top w_*(\cdot;x)$ decreases at rate $r_*$ indefinitely as $w_*(\cdot;x)$ keeps moving after hitting the origin.
When $w_*(T;x)\geq0$, let $x_\epsilon=x+(\epsilon,0,...,0)$ with $\epsilon>0$. Then
\begin{align*}
    \norm{w(T;x_\epsilon)-w(T;x)}_1&\geq\abs{\mathbf{1}^\top(w(T;x_\epsilon)- w(T;x))}\\
    &=\abs{\mathbf{1}^\top x_\epsilon-r_* T-\mathbf{1}^\top x+r_*T}\\
    &=\epsilon\\
    &=\norm{x_\epsilon-x}_1,
\end{align*}
so $Df(x)=1.$ Let $V_a(x)=\exp(a\mathbf{1}^\top x)$ where $a$ will be determined later. For $x\geq0,$
\begin{align*}
    KV_a(x)&=\E I(w_*(T;x)\geq0)V(w(T;x)+\bar{Z})\\
    &=\E I(w_*(T;x)\geq0)V(w_*(T;x)+\bar{Z})\\
    &\leq\E V(w_*(T;x)+\bar{Z})\\
    &=\E \exp(a\mathbf{1}^\top(w_*(T;x)+\bar{Z}))\\
    &=\E \exp(a(\mathbf{1}^\top x-r_*T+Z))\\
    &=V_a(x)\E\exp(a(Z-r_*T)),
\end{align*}
where the second equality is because $w(t;x)$ and $w_*(t;x)$ are the same until they hit the origin (boundary removal technique). Given $\E e^{\zeta Z}<\infty$, the optimal drift rate is
$$\lambda_*=\E\exp(a_*(Z-r_*T))=\inf_{a\in[0,\zeta]}\E\exp(a(Z-r_*T))<1.$$
By Theorem \ref{theorem_exponential},
\begin{align*}
    W_I(X_n,X_\infty)\leq &W_{V_{a_*}}(X_n,X_\infty)\\
    \leq&\frac{\lambda_*^n}{1-\lambda_*}\cdot\E d_{V_{a_*}}(X_0,X_1)\\
    \leq&\frac{\lambda_*^n}{1-\lambda_*}\cdot\E \int_0^1\exp\prs{a_*\mathbf{1}^\top\prs{(1-t)X_0+tX_1}}\norm{X_1-X_0}_1dt\\
    \leq&\frac{\lambda_*^n}{1-\lambda_*}\cdot\E\sbk{\norm{X_1-X_0}_1\int_0^1\exp\prs{a_*\prs{\mathbf{1}^\top X_0+t(\mathbf{1}^\top X_1-\mathbf{1}^\top X_0)}}dt}\\
    \leq&\frac{\lambda_*^n}{1-\lambda_*}\cdot\E\sbk{\norm{X_1-X_0}_1\frac{\exp\prs{a_*\prs{\mathbf{1}^\top X_0+t(\mathbf{1}^\top X_1-\mathbf{1}^\top X_0)}}}{a_*\prs{\mathbf{1}^\top X_1-\mathbf{1}^\top X_0}}\Bigg|_0^1}\\
    \leq&\frac{\lambda_*^n}{1-\lambda_*}\cdot\E\sbk{\norm{X_1-X_0}_1\frac{\exp\prs{a_*\mathbf{1}^\top X_1}-\exp\prs{a_*\mathbf{1}^\top X_0}}{a_*\mathbf{1}^\top X_1-a_*\mathbf{1}^\top X_0}},
\end{align*}
where the subscript $I$ corresponds to the intrinsic metric induced by $\norm{\cdot}_1$, which is $\norm{\cdot}_1$ itself. Since $\norm{\cdot}_1\geq\norm{\cdot}_2$, the above bound also holds for $W(X_n,X_\infty)$.
\end{proof}

%
%

\begin{acks}[Acknowledgments]
We would like to sincerely thank the anonymous referees for their insightful feedback, which has strengthened this paper. The material in this paper is partly supported by the Air Force Office of Scientific Research under award number FA9550-20-1-0397 and ONR N000142412655. Support from NSF 2229012, 2312204, 2403007 is also gratefully acknowledged.
\end{acks}
\bibliographystyle{imsart-nameyear} 
\bibliography{main}       

\begin{thebibliography}{73}

\bibitem[\protect\citeauthoryear{Albert and Chib}{1993}]{albert1993bayesian}
\begin{barticle}[author]
\bauthor{\bsnm{Albert},~\bfnm{James~H}\binits{J.~H.}} \AND \bauthor{\bsnm{Chib},~\bfnm{Siddhartha}\binits{S.}}
(\byear{1993}).
\btitle{Bayesian analysis of binary and polychotomous response data}.
\bjournal{Journal of the American statistical Association}
\bvolume{88}
\bpages{669--679}.
\end{barticle}
\endbibitem

\bibitem[\protect\citeauthoryear{Anari, Liu and Gharan}{2020}]{Anari20}
\begin{binproceedings}[author]
\bauthor{\bsnm{Anari},~\bfnm{Nima}\binits{N.}}, \bauthor{\bsnm{Liu},~\bfnm{Kuikui}\binits{K.}} \AND \bauthor{\bsnm{Gharan},~\bfnm{Shayan~Oveis}\binits{S.~O.}}
(\byear{2020}).
\btitle{Spectral independence in high-dimensional expanders and applications to the hardcore model}.
In \bbooktitle{2020 IEEE 61st Annual Symposium on Foundations of Computer Science (FOCS)}
\bpages{1319-1330}.
\end{binproceedings}
\endbibitem

\bibitem[\protect\citeauthoryear{Andrieu, Fort and Vihola}{2015}]{andrieu2015quantitative}
\begin{barticle}[author]
\bauthor{\bsnm{Andrieu},~\bfnm{Christophe}\binits{C.}}, \bauthor{\bsnm{Fort},~\bfnm{Gersende}\binits{G.}} \AND \bauthor{\bsnm{Vihola},~\bfnm{Matti}\binits{M.}}
(\byear{2015}).
\btitle{Quantitative convergence rates for subgeometric Markov chains}.
\bjournal{Journal of Applied Probability}
\bvolume{52}
\bpages{391--404}.
\end{barticle}
\endbibitem

\bibitem[\protect\citeauthoryear{Baxendale}{2005}]{baxendale2005renewal}
\begin{barticle}[author]
\bauthor{\bsnm{Baxendale},~\bfnm{Peter~H}\binits{P.~H.}}
(\byear{2005}).
\btitle{Renewal theory and computable convergence rates for geometrically ergodic Markov chains}.
\bjournal{The Annals of Applied Probability}
\bvolume{15}
\bpages{700--738}.
\end{barticle}
\endbibitem

\bibitem[\protect\citeauthoryear{Billingsley}{2013}]{billingsley2013convergence}
\begin{bbook}[author]
\bauthor{\bsnm{Billingsley},~\bfnm{Patrick}\binits{P.}}
(\byear{2013}).
\btitle{Convergence of Probability Measures}.
\bpublisher{John Wiley \& Sons}.
\end{bbook}
\endbibitem

\bibitem[\protect\citeauthoryear{Biswas, Jacob and Vanetti}{2019}]{biswas2019estimating}
\begin{barticle}[author]
\bauthor{\bsnm{Biswas},~\bfnm{Niloy}\binits{N.}}, \bauthor{\bsnm{Jacob},~\bfnm{Pierre~E}\binits{P.~E.}} \AND \bauthor{\bsnm{Vanetti},~\bfnm{Paul}\binits{P.}}
(\byear{2019}).
\btitle{Estimating convergence of {M}arkov chains with {L}-lag couplings}.
\bjournal{Advances in Neural Information Processing Systems}
\bvolume{32}.
\end{barticle}
\endbibitem

\bibitem[\protect\citeauthoryear{Bogachev}{2007}]{bogachev2007measure}
\begin{bbook}[author]
\bauthor{\bsnm{Bogachev},~\bfnm{Vladimir~I}\binits{V.~I.}}
(\byear{2007}).
\btitle{Measure Theory}
\bvolume{1}.
\bpublisher{Springer Science \& Business Media}.
\end{bbook}
\endbibitem

\bibitem[\protect\citeauthoryear{Bonferroni}{1936}]{bonferroni1936teoria}
\begin{barticle}[author]
\bauthor{\bsnm{Bonferroni},~\bfnm{Carlo}\binits{C.}}
(\byear{1936}).
\btitle{Teoria statistica delle classi e calcolo delle probabilita}.
\bjournal{Pubblicazioni del R Istituto Superiore di Scienze Economiche e Commericiali di Firenze}
\bvolume{8}
\bpages{3--62}.
\end{barticle}
\endbibitem

\bibitem[\protect\citeauthoryear{Bou-Rabee, Eberle and Zimmer}{2020}]{bou2020coupling}
\begin{barticle}[author]
\bauthor{\bsnm{Bou-Rabee},~\bfnm{Nawaf}\binits{N.}}, \bauthor{\bsnm{Eberle},~\bfnm{Andreas}\binits{A.}} \AND \bauthor{\bsnm{Zimmer},~\bfnm{Raphael}\binits{R.}}
(\byear{2020}).
\btitle{Coupling and convergence for {H}amiltonian {M}onte {C}arlo}.
\bjournal{The Annals of Applied Probability}
\bvolume{30}
\bpages{1209--1250}.
\end{barticle}
\endbibitem

\bibitem[\protect\citeauthoryear{Budhiraja and Lee}{2007}]{budhiraja2007long}
\begin{barticle}[author]
\bauthor{\bsnm{Budhiraja},~\bfnm{Amarjit}\binits{A.}} \AND \bauthor{\bsnm{Lee},~\bfnm{Chihoon}\binits{C.}}
(\byear{2007}).
\btitle{Long time asymptotics for constrained diffusions in polyhedral domains}.
\bjournal{Stochastic Processes and their Applications}
\bvolume{117}
\bpages{1014--1036}.
\end{barticle}
\endbibitem

\bibitem[\protect\citeauthoryear{Buraczewski, Damek and Mikosch}{2016}]{buraczewski2016stochastic}
\begin{bbook}[author]
\bauthor{\bsnm{Buraczewski},~\bfnm{Dariusz}\binits{D.}}, \bauthor{\bsnm{Damek},~\bfnm{Ewa}\binits{E.}} \AND \bauthor{\bsnm{Mikosch},~\bfnm{Thomas}\binits{T.}}
(\byear{2016}).
\btitle{Stochastic Models with Power-Law Tails: The Equation X= AX+ B}.
\bpublisher{Springer}.
\end{bbook}
\endbibitem

\bibitem[\protect\citeauthoryear{Burago et~al.}{2001}]{burago2001course}
\begin{bbook}[author]
\bauthor{\bsnm{Burago},~\bfnm{Dmitri}\binits{D.}}, \bauthor{\bsnm{Burago},~\bfnm{Yuri}\binits{Y.}}, \bauthor{\bsnm{Ivanov},~\bfnm{Sergei}\binits{S.}} \betal{et~al.}
(\byear{2001}).
\btitle{A course in metric geometry}
\bvolume{33}.
\bpublisher{American Mathematical Society Providence}.
\end{bbook}
\endbibitem

\bibitem[\protect\citeauthoryear{Butkovsky}{2014}]{butkovsky2014subgeometric}
\begin{barticle}[author]
\bauthor{\bsnm{Butkovsky},~\bfnm{Oleg}\binits{O.}}
(\byear{2014}).
\btitle{Subgeometric rates of convergence of Markov processes in the Wasserstein metric}.
\bjournal{The Annals of Applied Probability}
\bvolume{24}
\bpages{526--552}.
\end{barticle}
\endbibitem

\bibitem[\protect\citeauthoryear{Butkovsky, Kulik and Scheutzow}{2020}]{butkovsky2020generalized}
\begin{barticle}[author]
\bauthor{\bsnm{Butkovsky},~\bfnm{Oleg}\binits{O.}}, \bauthor{\bsnm{Kulik},~\bfnm{Alexei}\binits{A.}} \AND \bauthor{\bsnm{Scheutzow},~\bfnm{Michael}\binits{M.}}
(\byear{2020}).
\btitle{Generalized couplings and ergodic rates for {SPDE}s and other {M}arkov models}.
\bjournal{The Annals of Applied Probability}
\bvolume{30}
\bpages{1--39}.
\end{barticle}
\endbibitem

\bibitem[\protect\citeauthoryear{Chen and Eldan}{2022}]{Ronen22}
\begin{binproceedings}[author]
\bauthor{\bsnm{Chen},~\bfnm{Yuansi}\binits{Y.}} \AND \bauthor{\bsnm{Eldan},~\bfnm{Ronen}\binits{R.}}
(\byear{2022}).
\btitle{Localization schemes: A framework for proving mixing bounds for Markov chains (extended abstract)}.
In \bbooktitle{2022 IEEE 63rd Annual Symposium on Foundations of Computer Science (FOCS)}
\bpages{110-122}.
\end{binproceedings}
\endbibitem

\bibitem[\protect\citeauthoryear{Chen et~al.}{2001}]{chen2001fundamentals}
\begin{bbook}[author]
\bauthor{\bsnm{Chen},~\bfnm{Hong}\binits{H.}}, \bauthor{\bsnm{Yao},~\bfnm{David~D}\binits{D.~D.}} \betal{et~al.}
(\byear{2001}).
\btitle{Fundamentals of Queueing Networks: Performance, Asymptotics, and Optimization}
\bvolume{4}.
\bpublisher{Springer}.
\end{bbook}
\endbibitem

\bibitem[\protect\citeauthoryear{Diaconis and Freedman}{1999}]{diaconis1999iterated}
\begin{barticle}[author]
\bauthor{\bsnm{Diaconis},~\bfnm{Persi}\binits{P.}} \AND \bauthor{\bsnm{Freedman},~\bfnm{David}\binits{D.}}
(\byear{1999}).
\btitle{Iterated random functions}.
\bjournal{SIAM review}
\bvolume{41}
\bpages{45--76}.
\end{barticle}
\endbibitem

\bibitem[\protect\citeauthoryear{Dieuleveut, Durmus and Bach}{2020}]{dieuleveut2020bridging}
\begin{barticle}[author]
\bauthor{\bsnm{Dieuleveut},~\bfnm{Aymeric}\binits{A.}}, \bauthor{\bsnm{Durmus},~\bfnm{Alain}\binits{A.}} \AND \bauthor{\bsnm{Bach},~\bfnm{Francis}\binits{F.}}
(\byear{2020}).
\btitle{Bridging the gap between constant step size stochastic gradient descent and Markov chains}.
\bjournal{The Annals of Statistics}
\bvolume{48}
\bpages{1348--1382}.
\end{barticle}
\endbibitem

\bibitem[\protect\citeauthoryear{Douc et~al.}{2004}]{douc2004practical}
\begin{barticle}[author]
\bauthor{\bsnm{Douc},~\bfnm{Randal}\binits{R.}}, \bauthor{\bsnm{Fort},~\bfnm{Gersende}\binits{G.}}, \bauthor{\bsnm{Moulines},~\bfnm{Eric}\binits{E.}} \AND \bauthor{\bsnm{Soulier},~\bfnm{Philippe}\binits{P.}}
(\byear{2004}).
\btitle{Practical drift conditions for subgeometric rates of convergence}.
\bjournal{The Annals of Applied Probability}
\bvolume{14}
\bpages{1353--1377}.
\end{barticle}
\endbibitem

\bibitem[\protect\citeauthoryear{Douc et~al.}{2018}]{douc2018markov}
\begin{bbook}[author]
\bauthor{\bsnm{Douc},~\bfnm{Randal}\binits{R.}}, \bauthor{\bsnm{Moulines},~\bfnm{Eric}\binits{E.}}, \bauthor{\bsnm{Priouret},~\bfnm{Pierre}\binits{P.}} \AND \bauthor{\bsnm{Soulier},~\bfnm{Philippe}\binits{P.}}
(\byear{2018}).
\btitle{Markov Chains}.
\bpublisher{Springer}.
\end{bbook}
\endbibitem

\bibitem[\protect\citeauthoryear{Durmus, Fort and Moulines}{2016}]{durmus2016subgeometric}
\begin{barticle}[author]
\bauthor{\bsnm{Durmus},~\bfnm{Alain}\binits{A.}}, \bauthor{\bsnm{Fort},~\bfnm{Gersende}\binits{G.}} \AND \bauthor{\bsnm{Moulines},~\bfnm{{\'E}ric}\binits{{\'E}.}}
(\byear{2016}).
\btitle{{Subgeometric rates of convergence in Wasserstein distance for Markov chains}}.
\bjournal{Annales de l'Institut Henri Poincaré, Probabilités et Statistiques}
\bvolume{52}
\bpages{1799--1822}.
\end{barticle}
\endbibitem

\bibitem[\protect\citeauthoryear{Durmus and Moulines}{2015}]{durmus2015quantitative}
\begin{barticle}[author]
\bauthor{\bsnm{Durmus},~\bfnm{Alain}\binits{A.}} \AND \bauthor{\bsnm{Moulines},~\bfnm{{\'E}ric}\binits{{\'E}.}}
(\byear{2015}).
\btitle{Quantitative bounds of convergence for geometrically ergodic Markov chains in the Wasserstein distance with application to the Metropolis adjusted Langevin algorithm}.
\bjournal{Statistics and Computing}
\bvolume{25}
\bpages{5--19}.
\end{barticle}
\endbibitem

\bibitem[\protect\citeauthoryear{Durmus and Moulines}{2017}]{durmus2017nonasymptotic}
\begin{barticle}[author]
\bauthor{\bsnm{Durmus},~\bfnm{Alain}\binits{A.}} \AND \bauthor{\bsnm{Moulines},~\bfnm{{\'E}ric}\binits{{\'E}.}}
(\byear{2017}).
\btitle{Nonasymptotic convergence analysis for the unadjusted {L}angevin algorithm}.
\bjournal{The Annals of Applied Probability}
\bvolume{27}
\bpages{1551}.
\end{barticle}
\endbibitem

\bibitem[\protect\citeauthoryear{Durmus and Moulines}{2019}]{durmus2019high}
\begin{barticle}[author]
\bauthor{\bsnm{Durmus},~\bfnm{Alain}\binits{A.}} \AND \bauthor{\bsnm{Moulines},~\bfnm{{\'E}ric}\binits{{\'E}.}}
(\byear{2019}).
\btitle{High-dimensional {B}ayesian inference via the unadjusted {L}angevin algorithm}.
\bjournal{Bernoulli}
\bvolume{25}
\bpages{2854--2882}.
\end{barticle}
\endbibitem

\bibitem[\protect\citeauthoryear{Eberle}{2011}]{eberle2011reflection}
\begin{barticle}[author]
\bauthor{\bsnm{Eberle},~\bfnm{Andreas}\binits{A.}}
(\byear{2011}).
\btitle{Reflection coupling and {W}asserstein contractivity without convexity}.
\bjournal{Comptes Rendus Mathematique}
\bvolume{349}
\bpages{1101--1104}.
\end{barticle}
\endbibitem

\bibitem[\protect\citeauthoryear{Eberle}{2016}]{eberle2016reflection}
\begin{barticle}[author]
\bauthor{\bsnm{Eberle},~\bfnm{Andreas}\binits{A.}}
(\byear{2016}).
\btitle{Reflection couplings and contraction rates for diffusions}.
\bjournal{Probability Theory and Related Fields}
\bvolume{166}
\bpages{851--886}.
\end{barticle}
\endbibitem

\bibitem[\protect\citeauthoryear{Eberle, Guillin and Zimmer}{2019a}]{eberle2019couplings}
\begin{barticle}[author]
\bauthor{\bsnm{Eberle},~\bfnm{Andreas}\binits{A.}}, \bauthor{\bsnm{Guillin},~\bfnm{Arnaud}\binits{A.}} \AND \bauthor{\bsnm{Zimmer},~\bfnm{Raphael}\binits{R.}}
(\byear{2019}a).
\btitle{Couplings and quantitative contraction rates for {L}angevin dynamics}.
\bjournal{The Annals of Probability}
\bvolume{47}
\bpages{1982--2010}.
\end{barticle}
\endbibitem

\bibitem[\protect\citeauthoryear{Eberle, Guillin and Zimmer}{2019b}]{eberle2019quantitativeH}
\begin{barticle}[author]
\bauthor{\bsnm{Eberle},~\bfnm{Andreas}\binits{A.}}, \bauthor{\bsnm{Guillin},~\bfnm{Arnaud}\binits{A.}} \AND \bauthor{\bsnm{Zimmer},~\bfnm{Raphael}\binits{R.}}
(\byear{2019}b).
\btitle{Quantitative {H}arris-type theorems for diffusions and {M}c{K}ean--{V}lasov processes}.
\bjournal{Transactions of the American Mathematical Society}
\bvolume{371}
\bpages{7135--7173}.
\end{barticle}
\endbibitem

\bibitem[\protect\citeauthoryear{Eberle and Majka}{2019}]{eberle2019quantitative}
\begin{barticle}[author]
\bauthor{\bsnm{Eberle},~\bfnm{Andreas}\binits{A.}} \AND \bauthor{\bsnm{Majka},~\bfnm{Mateusz~B}\binits{M.~B.}}
(\byear{2019}).
\btitle{Quantitative contraction rates for {M}arkov chains on general state spaces}.
\bjournal{Electronic Journal of Probability}
\bvolume{24}
\bpages{26}.
\end{barticle}
\endbibitem

\bibitem[\protect\citeauthoryear{Eberle and Zimmer}{2019}]{eberle2019sticky}
\begin{binproceedings}[author]
\bauthor{\bsnm{Eberle},~\bfnm{Andreas}\binits{A.}} \AND \bauthor{\bsnm{Zimmer},~\bfnm{Raphael}\binits{R.}}
(\byear{2019}).
\btitle{Sticky couplings of multidimensional diffusions with different drifts}.
In \bbooktitle{Annales de l’Institut Henri Poincar{\'e}-Probabilit{\'e}s et Statistiques}
\bvolume{55}
\bpages{2370--2394}.
\end{binproceedings}
\endbibitem

\bibitem[\protect\citeauthoryear{Ghosh and Marecek}{2022}]{ghosh2022iterated}
\begin{barticle}[author]
\bauthor{\bsnm{Ghosh},~\bfnm{Ramen}\binits{R.}} \AND \bauthor{\bsnm{Marecek},~\bfnm{Jakub}\binits{J.}}
(\byear{2022}).
\btitle{Iterated function systems: A comprehensive survey}.
\bjournal{arXiv preprint arXiv:2211.14661}.
\end{barticle}
\endbibitem

\bibitem[\protect\citeauthoryear{Gibbs}{2004}]{gibbs2004convergence}
\begin{barticle}[author]
\bauthor{\bsnm{Gibbs},~\bfnm{Alison~L}\binits{A.~L.}}
(\byear{2004}).
\btitle{Convergence in the Wasserstein metric for Markov chain Monte Carlo algorithms with applications to image restoration}.
\bjournal{Stochastic Models}
\bvolume{20}
\bpages{473--492}.
\end{barticle}
\endbibitem

\bibitem[\protect\citeauthoryear{Glynn and Wang}{2018}]{glynn2018rate}
\begin{barticle}[author]
\bauthor{\bsnm{Glynn},~\bfnm{Peter~W}\binits{P.~W.}} \AND \bauthor{\bsnm{Wang},~\bfnm{Rob~J}\binits{R.~J.}}
(\byear{2018}).
\btitle{On the rate of convergence to equilibrium for reflected Brownian motion}.
\bjournal{Queueing Systems}
\bvolume{89}
\bpages{165--197}.
\end{barticle}
\endbibitem

\bibitem[\protect\citeauthoryear{Gould}{1972}]{gould1972combinatorial}
\begin{bbook}[author]
\bauthor{\bsnm{Gould},~\bfnm{H.~W.}\binits{H.~W.}}
(\byear{1972}).
\btitle{Combinatorial Identities: A Standardized Set of Tables Listing 500 Binomial Coefficient Summations}.
\bpublisher{Morgantown Printing and Binding Co.}
\end{bbook}
\endbibitem

\bibitem[\protect\citeauthoryear{Hairer and Mattingly}{2008}]{hairer2008spectral}
\begin{barticle}[author]
\bauthor{\bsnm{Hairer},~\bfnm{Martin}\binits{M.}} \AND \bauthor{\bsnm{Mattingly},~\bfnm{Jonathan~C}\binits{J.~C.}}
(\byear{2008}).
\btitle{Spectral gaps in {W}asserstein distances and the 2{D} stochastic {N}avier--{S}tokes equations}.
\bjournal{The Annals of Probability}
\bvolume{36}
\bpages{2050--2091}.
\end{barticle}
\endbibitem

\bibitem[\protect\citeauthoryear{Hairer and Mattingly}{2011}]{hairer2011yet}
\begin{binproceedings}[author]
\bauthor{\bsnm{Hairer},~\bfnm{Martin}\binits{M.}} \AND \bauthor{\bsnm{Mattingly},~\bfnm{Jonathan~C}\binits{J.~C.}}
(\byear{2011}).
\btitle{Yet another look at {H}arris’ ergodic theorem for {M}arkov chains}.
In \bbooktitle{Seminar on Stochastic Analysis, Random Fields and Applications VI: Centro Stefano Franscini, Ascona, May 2008}
\bpages{109--117}.
\bpublisher{Springer}.
\end{binproceedings}
\endbibitem

\bibitem[\protect\citeauthoryear{Hairer, Mattingly and Scheutzow}{2011}]{hairer2011asymptotic}
\begin{barticle}[author]
\bauthor{\bsnm{Hairer},~\bfnm{Martin}\binits{M.}}, \bauthor{\bsnm{Mattingly},~\bfnm{Jonathan~C}\binits{J.~C.}} \AND \bauthor{\bsnm{Scheutzow},~\bfnm{Michael}\binits{M.}}
(\byear{2011}).
\btitle{Asymptotic coupling and a general form of Harris’ theorem with applications to stochastic delay equations}.
\bjournal{Probability Theory and Related Fields}
\bvolume{149}
\bpages{223--259}.
\end{barticle}
\endbibitem

\bibitem[\protect\citeauthoryear{Hairer, Stuart and Vollmer}{2014}]{hairer2014spectral}
\begin{barticle}[author]
\bauthor{\bsnm{Hairer},~\bfnm{Martin}\binits{M.}}, \bauthor{\bsnm{Stuart},~\bfnm{AM}\binits{A.}} \AND \bauthor{\bsnm{Vollmer},~\bfnm{Sebastian}\binits{S.}}
(\byear{2014}).
\btitle{Spectral gaps for a {M}etropolis--{H}astings algorithm in infinite dimensions}.
\bjournal{The Annals of Applied Probability}
\bvolume{24}
\bpages{2455--2490}.
\end{barticle}
\endbibitem

\bibitem[\protect\citeauthoryear{Harrison and Reiman}{1981}]{harrison1981reflected}
\begin{barticle}[author]
\bauthor{\bsnm{Harrison},~\bfnm{J~Michael}\binits{J.~M.}} \AND \bauthor{\bsnm{Reiman},~\bfnm{Martin~I}\binits{M.~I.}}
(\byear{1981}).
\btitle{Reflected {B}rownian motion on an orthant}.
\bjournal{The Annals of Probability}
\bvolume{9}
\bpages{302--308}.
\end{barticle}
\endbibitem

\bibitem[\protect\citeauthoryear{Hu and Kirk}{1978}]{hu1978local}
\begin{barticle}[author]
\bauthor{\bsnm{Hu},~\bfnm{Thakyin}\binits{T.}} \AND \bauthor{\bsnm{Kirk},~\bfnm{WA}\binits{W.}}
(\byear{1978}).
\btitle{Local contractions in metric spaces}.
\bjournal{Proceedings of the American Mathematical Society}
\bvolume{68}
\bpages{121--124}.
\end{barticle}
\endbibitem

\bibitem[\protect\citeauthoryear{Jarner and Roberts}{2002}]{jarner2002polynomial}
\begin{barticle}[author]
\bauthor{\bsnm{Jarner},~\bfnm{S{\o}ren~F}\binits{S.~F.}} \AND \bauthor{\bsnm{Roberts},~\bfnm{Gareth~O}\binits{G.~O.}}
(\byear{2002}).
\btitle{Polynomial convergence rates of Markov chains}.
\bjournal{The Annals of Applied Probability}
\bvolume{12}
\bpages{224--247}.
\end{barticle}
\endbibitem

\bibitem[\protect\citeauthoryear{Jarner and Tweedie}{2001}]{jarner2001locally}
\begin{barticle}[author]
\bauthor{\bsnm{Jarner},~\bfnm{SF}\binits{S.}} \AND \bauthor{\bsnm{Tweedie},~\bfnm{RL}\binits{R.}}
(\byear{2001}).
\btitle{Locally contracting iterated functions and stability of {M}arkov chains}.
\bjournal{Journal of Applied Probability}
\bvolume{38}
\bpages{494--507}.
\end{barticle}
\endbibitem

\bibitem[\protect\citeauthoryear{Jones and Hobert}{2001}]{jones2001honest}
\begin{barticle}[author]
\bauthor{\bsnm{Jones},~\bfnm{Galin~L.}\binits{G.~L.}} \AND \bauthor{\bsnm{Hobert},~\bfnm{James~P.}\binits{J.~P.}}
(\byear{2001}).
\btitle{Honest exploration of intractable probability distributions via {M}arkov chain {M}onte {C}arlo}.
\bjournal{Statistical Science}
\bvolume{16}
\bpages{312--334}.
\end{barticle}
\endbibitem

\bibitem[\protect\citeauthoryear{Kella and Whitt}{1992}]{kella1992tandem}
\begin{barticle}[author]
\bauthor{\bsnm{Kella},~\bfnm{Offer}\binits{O.}} \AND \bauthor{\bsnm{Whitt},~\bfnm{Ward}\binits{W.}}
(\byear{1992}).
\btitle{A tandem fluid network with {L}\'evy input}.
\bjournal{Queueing and Related Models}
\bpages{112--128}.
\end{barticle}
\endbibitem

\bibitem[\protect\citeauthoryear{Kiefer and Wolfowitz}{1956}]{kiefer1956characteristics}
\begin{barticle}[author]
\bauthor{\bsnm{Kiefer},~\bfnm{J}\binits{J.}} \AND \bauthor{\bsnm{Wolfowitz},~\bfnm{J}\binits{J.}}
(\byear{1956}).
\btitle{On the characteristics of the general queueing process, with applications to random walk}.
\bjournal{The Annals of Mathematical Statistics}
\bpages{147--161}.
\end{barticle}
\endbibitem

\bibitem[\protect\citeauthoryear{Lazi and Sandri}{2021}]{lazi2021sub}
\begin{barticle}[author]
\bauthor{\bsnm{Lazi},~\bfnm{Petra}\binits{P.}} \AND \bauthor{\bsnm{Sandri},~\bfnm{Nikola}\binits{N.}}
(\byear{2021}).
\btitle{On sub-geometric ergodicity of diffusion processes}.
\bjournal{Bernoulli}
\bvolume{27}
\bpages{348--380}.
\end{barticle}
\endbibitem

\bibitem[\protect\citeauthoryear{Madras and Sezer}{2010}]{madras2010quantitative}
\begin{barticle}[author]
\bauthor{\bsnm{Madras},~\bfnm{Neal}\binits{N.}} \AND \bauthor{\bsnm{Sezer},~\bfnm{Deniz}\binits{D.}}
(\byear{2010}).
\btitle{Quantitative bounds for Markov chain convergence: Wasserstein and total variation distances}.
\bjournal{Bernoulli}
\bvolume{16}
\bpages{882--908}.
\end{barticle}
\endbibitem

\bibitem[\protect\citeauthoryear{Mangoubi and Smith}{2019}]{mangoubi2019mixing}
\begin{binproceedings}[author]
\bauthor{\bsnm{Mangoubi},~\bfnm{Oren}\binits{O.}} \AND \bauthor{\bsnm{Smith},~\bfnm{Aaron}\binits{A.}}
(\byear{2019}).
\btitle{Mixing of {H}amiltonian {M}onte {C}arlo on strongly log-concave distributions 2: Numerical integrators}.
In \bbooktitle{Proceedings of the Twenty-Second International Conference on Artificial Intelligence and Statistics}
\bvolume{89}
\bpages{586--595}.
\bpublisher{PMLR}.
\end{binproceedings}
\endbibitem

\bibitem[\protect\citeauthoryear{Meyn and Tweedie}{1994}]{meyn1994computable}
\begin{barticle}[author]
\bauthor{\bsnm{Meyn},~\bfnm{Sean~P}\binits{S.~P.}} \AND \bauthor{\bsnm{Tweedie},~\bfnm{Robert~L}\binits{R.~L.}}
(\byear{1994}).
\btitle{Computable bounds for geometric convergence rates of Markov chains}.
\bjournal{The Annals of Applied Probability}
\bvolume{4}
\bpages{981--1011}.
\end{barticle}
\endbibitem

\bibitem[\protect\citeauthoryear{Meyn and Tweedie}{2009}]{meyn_tweedie_glynn_2009}
\begin{bbook}[author]
\bauthor{\bsnm{Meyn},~\bfnm{Sean~P}\binits{S.~P.}} \AND \bauthor{\bsnm{Tweedie},~\bfnm{Richard~L}\binits{R.~L.}}
(\byear{2009}).
\btitle{Markov Chains and Stochastic Stability},
\bedition{2} ed.
\bseries{Cambridge Mathematical Library}.
\bpublisher{Cambridge University Press}.
\end{bbook}
\endbibitem

\bibitem[\protect\citeauthoryear{Monmarch{\'e}}{2021}]{monmarche2021high}
\begin{barticle}[author]
\bauthor{\bsnm{Monmarch{\'e}},~\bfnm{Pierre}\binits{P.}}
(\byear{2021}).
\btitle{High-dimensional {MCMC} with a standard splitting scheme for the underdamped {L}angevin diffusion.}
\bjournal{Electronic Journal of Statistics}
\bvolume{15}
\bpages{4117--4166}.
\end{barticle}
\endbibitem

\bibitem[\protect\citeauthoryear{Nguyen}{2024}]{nguyen2024polynomial}
\begin{barticle}[author]
\bauthor{\bsnm{Nguyen},~\bfnm{Hung~D}\binits{H.~D.}}
(\byear{2024}).
\btitle{Polynomial mixing of a stochastic wave equation with dissipative damping}.
\bjournal{Applied Mathematics \& Optimization}
\bvolume{89}
\bpages{21}.
\end{barticle}
\endbibitem

\bibitem[\protect\citeauthoryear{Ollivier}{2009}]{ollivier2009ricci}
\begin{barticle}[author]
\bauthor{\bsnm{Ollivier},~\bfnm{Yann}\binits{Y.}}
(\byear{2009}).
\btitle{Ricci curvature of {M}arkov chains on metric spaces}.
\bjournal{Journal of Functional Analysis}
\bvolume{256}
\bpages{810--864}.
\end{barticle}
\endbibitem

\bibitem[\protect\citeauthoryear{Qin and Hobert}{2021}]{qin2021limitations}
\begin{barticle}[author]
\bauthor{\bsnm{Qin},~\bfnm{Qian}\binits{Q.}} \AND \bauthor{\bsnm{Hobert},~\bfnm{James~P}\binits{J.~P.}}
(\byear{2021}).
\btitle{On the limitations of single-step drift and minorization in Markov chain convergence analysis}.
\bjournal{The Annals of Applied Probability}
\bvolume{31}
\bpages{1633--1659}.
\end{barticle}
\endbibitem

\bibitem[\protect\citeauthoryear{Qin and Hobert}{2022a}]{qin2022wasserstein}
\begin{barticle}[author]
\bauthor{\bsnm{Qin},~\bfnm{Qian}\binits{Q.}} \AND \bauthor{\bsnm{Hobert},~\bfnm{James~P}\binits{J.~P.}}
(\byear{2022}a).
\btitle{Wasserstein-based methods for convergence complexity analysis of MCMC with applications}.
\bjournal{The Annals of Applied Probability}
\bvolume{32}
\bpages{124--166}.
\end{barticle}
\endbibitem

\bibitem[\protect\citeauthoryear{Qin and Hobert}{2022b}]{qin2022geometric}
\begin{barticle}[author]
\bauthor{\bsnm{Qin},~\bfnm{Qian}\binits{Q.}} \AND \bauthor{\bsnm{Hobert},~\bfnm{James~P.}\binits{J.~P.}}
(\byear{2022}b).
\btitle{{Geometric convergence bounds for Markov chains in Wasserstein distance based on generalized drift and contraction conditions}}.
\bjournal{Annales de l'Institut Henri Poincaré, Probabilités et Statistiques}
\bvolume{58}
\bpages{872--889}.
\end{barticle}
\endbibitem

\bibitem[\protect\citeauthoryear{Qu, Blanchet and Glynn}{2024}]{qu2024deep}
\begin{barticle}[author]
\bauthor{\bsnm{Qu},~\bfnm{Yanlin}\binits{Y.}}, \bauthor{\bsnm{Blanchet},~\bfnm{Jose}\binits{J.}} \AND \bauthor{\bsnm{Glynn},~\bfnm{Peter~W}\binits{P.~W.}}
(\byear{2024}).
\btitle{Deep learning for computing convergence rates of Markov chains}.
\bjournal{Advances in Neural Information Processing Systems}
\bvolume{37}
\bpages{84777--84798}.
\end{barticle}
\endbibitem

\bibitem[\protect\citeauthoryear{Raissi, Perdikaris and Karniadakis}{2019}]{raissi2019physics}
\begin{barticle}[author]
\bauthor{\bsnm{Raissi},~\bfnm{Maziar}\binits{M.}}, \bauthor{\bsnm{Perdikaris},~\bfnm{Paris}\binits{P.}} \AND \bauthor{\bsnm{Karniadakis},~\bfnm{George~E}\binits{G.~E.}}
(\byear{2019}).
\btitle{Physics-informed neural networks: A deep learning framework for solving forward and inverse problems involving nonlinear partial differential equations}.
\bjournal{Journal of Computational Physics}
\bvolume{378}
\bpages{686--707}.
\end{barticle}
\endbibitem

\bibitem[\protect\citeauthoryear{Rajaratnam and Sparks}{2015}]{rajaratnam2015mcmc}
\begin{barticle}[author]
\bauthor{\bsnm{Rajaratnam},~\bfnm{Bala}\binits{B.}} \AND \bauthor{\bsnm{Sparks},~\bfnm{Doug}\binits{D.}}
(\byear{2015}).
\btitle{{MCMC}-based inference in the era of big data: A fundamental analysis of the convergence complexity of high-dimensional chains}.
\bjournal{arXiv preprint arXiv:1508.00947}.
\end{barticle}
\endbibitem

\bibitem[\protect\citeauthoryear{Rosenthal}{1995}]{rosenthal1995minorization}
\begin{barticle}[author]
\bauthor{\bsnm{Rosenthal},~\bfnm{Jeffrey~S}\binits{J.~S.}}
(\byear{1995}).
\btitle{Minorization conditions and convergence rates for Markov chain Monte Carlo}.
\bjournal{Journal of the American Statistical Association}
\bvolume{90}
\bpages{558--566}.
\end{barticle}
\endbibitem

\bibitem[\protect\citeauthoryear{Sandri{\'c}, Arapostathis and Pang}{2022}]{sandric2022subexponential}
\begin{barticle}[author]
\bauthor{\bsnm{Sandri{\'c}},~\bfnm{Nikola}\binits{N.}}, \bauthor{\bsnm{Arapostathis},~\bfnm{Ari}\binits{A.}} \AND \bauthor{\bsnm{Pang},~\bfnm{Guodong}\binits{G.}}
(\byear{2022}).
\btitle{Subexponential upper and lower bounds in {W}asserstein distance for {M}arkov processes}.
\bjournal{Applied Mathematics \& Optimization}
\bvolume{85}
\bpages{37}.
\end{barticle}
\endbibitem

\bibitem[\protect\citeauthoryear{Simsekli, Sagun and Gurbuzbalaban}{2019}]{simsekli2019tail}
\begin{binproceedings}[author]
\bauthor{\bsnm{Simsekli},~\bfnm{Umut}\binits{U.}}, \bauthor{\bsnm{Sagun},~\bfnm{Levent}\binits{L.}} \AND \bauthor{\bsnm{Gurbuzbalaban},~\bfnm{Mert}\binits{M.}}
(\byear{2019}).
\btitle{A tail-index analysis of stochastic gradient noise in deep neural networks}.
In \bbooktitle{International Conference on Machine Learning}
\bpages{5827--5837}.
\bpublisher{PMLR}.
\end{binproceedings}
\endbibitem

\bibitem[\protect\citeauthoryear{Sirignano and Spiliopoulos}{2018}]{sirignano2018dgm}
\begin{barticle}[author]
\bauthor{\bsnm{Sirignano},~\bfnm{Justin}\binits{J.}} \AND \bauthor{\bsnm{Spiliopoulos},~\bfnm{Konstantinos}\binits{K.}}
(\byear{2018}).
\btitle{{DGM}: A deep learning algorithm for solving partial differential equations}.
\bjournal{Journal of Computational Physics}
\bvolume{375}
\bpages{1339--1364}.
\end{barticle}
\endbibitem

\bibitem[\protect\citeauthoryear{Spitzer}{1956}]{spitzer1956combinatorial}
\begin{barticle}[author]
\bauthor{\bsnm{Spitzer},~\bfnm{Frank}\binits{F.}}
(\byear{1956}).
\btitle{A combinatorial lemma and its application to probability theory}.
\bjournal{Transactions of the American Mathematical Society}
\bvolume{82}
\bpages{323--339}.
\end{barticle}
\endbibitem

\bibitem[\protect\citeauthoryear{Stein and Shakarchi}{2009}]{stein2009real}
\begin{bbook}[author]
\bauthor{\bsnm{Stein},~\bfnm{Elias~M}\binits{E.~M.}} \AND \bauthor{\bsnm{Shakarchi},~\bfnm{Rami}\binits{R.}}
(\byear{2009}).
\btitle{Real Analysis: Measure Theory, Integration, and Hilbert Spaces}.
\bpublisher{Princeton University Press}.
\end{bbook}
\endbibitem

\bibitem[\protect\citeauthoryear{Steinsaltz}{1999}]{steinsaltz1999locally}
\begin{barticle}[author]
\bauthor{\bsnm{Steinsaltz},~\bfnm{David}\binits{D.}}
(\byear{1999}).
\btitle{Locally contractive iterated function systems}.
\bjournal{The Annals of Probability}
\bpages{1952--1979}.
\end{barticle}
\endbibitem

\bibitem[\protect\citeauthoryear{Stenflo}{2001}]{stenflo2001markov}
\begin{barticle}[author]
\bauthor{\bsnm{Stenflo},~\bfnm{{\"O}rjan}\binits{{\"O}.}}
(\byear{2001}).
\btitle{Markov chains in random environments and random iterated function systems}.
\bjournal{Transactions of the American Mathematical Society}
\bvolume{353}
\bpages{3547--3562}.
\end{barticle}
\endbibitem

\bibitem[\protect\citeauthoryear{Stenflo}{2012}]{stenflo2012survey}
\begin{barticle}[author]
\bauthor{\bsnm{Stenflo},~\bfnm{{\"O}rjan}\binits{{\"O}.}}
(\byear{2012}).
\btitle{A survey of average contractive iterated function systems}.
\bjournal{Journal of Difference Equations and Applications}
\bvolume{18}
\bpages{1355--1380}.
\end{barticle}
\endbibitem

\bibitem[\protect\citeauthoryear{Tuominen and Tweedie}{1994}]{tuominen1994subgeometric}
\begin{barticle}[author]
\bauthor{\bsnm{Tuominen},~\bfnm{Pekka}\binits{P.}} \AND \bauthor{\bsnm{Tweedie},~\bfnm{Richard~L}\binits{R.~L.}}
(\byear{1994}).
\btitle{Subgeometric rates of convergence of f-ergodic Markov chains}.
\bjournal{Advances in Applied Probability}
\bvolume{26}
\bpages{775--798}.
\end{barticle}
\endbibitem

\bibitem[\protect\citeauthoryear{Villani et~al.}{2009}]{villani2009optimal}
\begin{bbook}[author]
\bauthor{\bsnm{Villani},~\bfnm{C{\'e}dric}\binits{C.}} \betal{et~al.}
(\byear{2009}).
\btitle{Optimal Transport: Old and New}
\bvolume{338}.
\bpublisher{Springer}.
\end{bbook}
\endbibitem

\bibitem[\protect\citeauthoryear{Yu et~al.}{2021}]{yu2021analysis}
\begin{barticle}[author]
\bauthor{\bsnm{Yu},~\bfnm{Lu}\binits{L.}}, \bauthor{\bsnm{Balasubramanian},~\bfnm{Krishnakumar}\binits{K.}}, \bauthor{\bsnm{Volgushev},~\bfnm{Stanislav}\binits{S.}} \AND \bauthor{\bsnm{Erdogdu},~\bfnm{Murat~A}\binits{M.~A.}}
(\byear{2021}).
\btitle{An analysis of constant step size SGD in the non-convex regime: Asymptotic normality and bias}.
\bjournal{Advances in Neural Information Processing Systems}
\bvolume{34}
\bpages{4234--4248}.
\end{barticle}
\endbibitem

\bibitem[\protect\citeauthoryear{Zhou et~al.}{2022}]{zhou2022dimension}
\begin{barticle}[author]
\bauthor{\bsnm{Zhou},~\bfnm{Quan}\binits{Q.}}, \bauthor{\bsnm{Yang},~\bfnm{Jun}\binits{J.}}, \bauthor{\bsnm{Vats},~\bfnm{Dootika}\binits{D.}}, \bauthor{\bsnm{Roberts},~\bfnm{Gareth~O}\binits{G.~O.}} \AND \bauthor{\bsnm{Rosenthal},~\bfnm{Jeffrey~S}\binits{J.~S.}}
(\byear{2022}).
\btitle{Dimension-free mixing for high-dimensional Bayesian variable selection}.
\bjournal{Journal of the Royal Statistical Society: Series B (Statistical Methodology)}
\bvolume{84}
\bpages{1751--1784}.
\end{barticle}
\endbibitem

\bibitem[\protect\citeauthoryear{Zimmer}{2017}]{zimmer2017explicit}
\begin{barticle}[author]
\bauthor{\bsnm{Zimmer},~\bfnm{Raphael}\binits{R.}}
(\byear{2017}).
\btitle{Explicit contraction rates for a class of degenerate and infinite-dimensional diffusions}.
\bjournal{Stochastics and Partial Differential Equations: Analysis and Computations}
\bvolume{5}
\bpages{368--399}.
\end{barticle}
\endbibitem

\end{thebibliography}


\end{document}